\documentclass[12pt, reqno]{amsart}
\usepackage{amsmath, amstext, amsbsy, amssymb, amscd}
\usepackage{amsmath}
\usepackage{amsxtra}
\usepackage{amscd}
\usepackage{amsthm}
\usepackage{amsfonts}
\usepackage{amssymb}
\usepackage{eucal}
\usepackage{color}

\newtheorem{theorem}{Theorem}[section]

\newtheorem{lemma}[theorem]{Lemma}
\newtheorem{proposition}[theorem]{Proposition}
\newtheorem{corollary}[theorem]{Corollary}
\theoremstyle{definition}
\newtheorem{definition}[theorem]{Definition}

\theoremstyle{remark}
\newtheorem{remark}[theorem]{Remark}

\definecolor{A}{rgb}{.75,1,.75}

\numberwithin{equation}{section}

\newcommand{\affH}{\widehat{{\mathcal H}}_n}
\newcommand{\C}{ \mathbb C }
\newcommand{\Cl}{ {\mathcal C} }
\newcommand{\ep}{\varepsilon}
\newcommand{\g}{\gamma}
\newcommand{\hf}{\frac12}
\newcommand{\HCaff}{\widehat{\mathcal H}c_n}
\newcommand{\HCn}{{{\mathcal H}c}_n}
\newcommand{\Hn}{{\mathcal H}_n^-}
\newcommand{\Hcycl}{{\mathcal H}_{n}^{I,-}}
\newcommand{\Haff}{\widehat{{\mathcal H}}_{n}^{-}}
\newcommand{\tHn}{{\mathcal H}^\thicksim_n}
\newcommand{\tHaff}{\widehat{{\mathcal H}}_n^\thicksim}
\newcommand{\tS}{S^\thicksim}
\newcommand{\tT}{\tilde{T}}

\newcommand{\pp}{p}
\newcommand{\qq}{q}
\newcommand{\PP}{P}
\newcommand{\QQ}{Q}
\newcommand{\rr}{R}
\newcommand{\tprime}{{\mathcal T}}

\newcommand{\Z}{ \mathbb Z }

\newcommand{\An}{ \mathcal{A} }
\newcommand{\Bn}{ \mathcal{B} }

{\vskip-\lastskip\medskip
  \noindent
  {\em #1.}\enspace
  }%
{\qed\par\medskip
  }

\begin{document}
\title{Spin Hecke algebras of finite and affine types}

\author[Weiqiang Wang]{Weiqiang Wang}
\address{Department of Mathematics, University of Virginia,
Charlottesville, VA 22904} \email{ww9c@virginia.edu}



\begin{abstract}
We introduce the spin Hecke algebra, which is a $q$-deformation of
the spin symmetric group algebra, and its affine generalization.
We establish an algebra isomorphism which relates our spin
(affine) Hecke algebras to the (affine) Hecke-Clifford algebras of
Olshanski and Jones-Nazarov. Relation between the spin (affine)
Hecke algebra and a nonstandard presentation of the usual (affine)
Hecke algebra is displayed, and the notion of covering (affine)
Hecke algebra is introduced to provide a link between these
algebras. Various algebraic structures for the spin (affine) Hecke
algebra are established.
\end{abstract}

\maketitle

\date{}

\section{Introduction}
\subsection{A basic question}
The spin (or projective) representations of the symmetric group
were first developed by I.~Schur \cite{Sch} in 1911. We refer to
J\'ozefiak \cite{Joz} for an excellent modern exposition of
Schur's work by a systematic use of superalgebras. The symmetric
group $S_n$ admits a double cover $\tS_n$, nontrivial for $n \ge
4$:
\begin{equation}  \label{centralext}
1 \longrightarrow\mathbb Z_2  {\longrightarrow} \tS_n
 {\longrightarrow} S_n \longrightarrow 1.
\end{equation}
A spin representation of $S_n$ is equivalent to a representation
of the algebra $\C S_n^- :=\C \tS_n/\langle z+1 \rangle$, the
quotient of the group algebra $\C\tS_n$ by the ideal $\langle z+1
\rangle$, where $z$ denotes the central generator of order $2$
coming from $\Z_2$. The algebra $\C S_n^-$ has a presentation with
generators $t_i $ $(1\le i \le n-1)$ subject to the relations:
\begin{eqnarray}
 t_i^2=1, && t_it_{i+1}t_i =t_{i+1}t_it_{i+1},  \label{eq:braidt}
\\
 t_it_j=-t_jt_i  &&   (|i -j|>1). \label{eq:tcomm}
\end{eqnarray}

As is well known, Hecke algebras have played an important role in
various aspects of representation theory (for some recent
developments see the books of Ariki \cite{Ar} and Kleshchev
\cite{Kle} and the references therein). We ask the following basic
question: is there a natural $q$-deformation (i.e. Hecke algebra)
for $\C S_n^-$ and $\C \tS_n$? It is conceivable that a canonical
solution to this question, if it exists, might open the door to
further new developments in representation theory.

However, there is no standard procedure to define Hecke algebras
except for Coxeter groups and perhaps for complex reflection
groups. The group $\tS_n$ is neither a Coxeter group nor a complex
reflection group.
\subsection{An affirmative answer}
In this paper we introduce the spin Hecke algebra $\Hn$ and the
covering Hecke algebra $\tHn$ as $q$-deformations of $\C S_n^-$ and
$\C\tS_n$ respectively. We also introduce the spin and the covering
affine Hecke algebras, denoted by $\Haff$ and $\tHaff$. The spin
(affine) Hecke algebras arise from different setups and they enjoy
various favorable properties.

Set
$$\ep =q -q^{-1}.$$
The spin affine Hecke algebra $\Haff$ is the $\C(q)$-algebra
generated by $\rr_i$, $1 \le i \le n-1$, and $p_i, q_i$, $1 \le i
\le n$, subject to the following relations:
\begin{eqnarray*}
\rr_i^2 &=& -\ep^2-2    \\
\rr_i \rr_j &=& -\rr_j \rr_i  \quad (|i-j|>1)   \\
\rr_i \rr_{i+1} \rr_i - \rr_{i+1} \rr_i \rr_{i+1} &=& \ep^2
(\rr_{i+1} -\rr_{i}) \\
\rr_i \pp_i &=& \pp_{i+1}\rr_i +\ep (\qq_i -\qq_{i+1})   \\
\rr_i \qq_i &=& -\qq_{i+1}\rr_i -\ep (\pp_i +\pp_{i+1})\\
%
\rr_i \pp_{j} = \pp_{j}\rr_i, \;\;\;
 \rr_i \qq_{j} &=& -\qq_{j}\rr_i \quad (j\neq i,i+1) \\
\pp_i \pp_j = \pp_j\pp_i,
 \;\;\;\; \qq_i \qq_j &=& -\qq_j\qq_i  \quad\; (i \neq j)   \\
\pp_i^2 +\qq_i^2 &=&1 \\
\pp_i \qq_j &=& \qq_j\pp_i  \quad (\forall i, j).
\end{eqnarray*}

The algebra $\Haff$ can be viewed as a quantum version of the
degenerate spin affine Hecke algebra introduced in \cite{Wa} (see
Section~\ref{subsec:degen}). The subalgebra generated by $\rr_i$
$(1\le i \le n-1)$ is the spin Hecke algebra $\Hn$ of finite type.
Among the noteworthy features of $\Hn$ and $\Haff$ are the deformed
braid relations and the two dependent sets of loop generators. Both
$\Haff$ and $\Hn$ admit superalgebra structures with each $\pp_i$
being even and $\qq_i, \rr_i$ being odd.

\subsection{Several related algebras}
To formulate a certain Schur-Jimbo type duality,
Olshanski \cite{Ol} introduced a Hecke-Clifford algebra $\HCn$,
which is a $q$-deformation of the semidirect product $\Cl_n
\rtimes \C S_n$ and is generated by the usual Hecke algebra
$\mathcal H_n$ for $S_n$ and the Clifford algebra $\Cl_n$ in $n$
variables. The affine Hecke-Clifford algebra $\HCaff$ was
introduced by Jones-Nazarov \cite{JN} to study the $q$-Young
symmetrizer for $\HCn$, and the modular representation theory of
$\HCaff$ has been developed by Brundan-Kleshchev \cite{BK}. A
degenerate version of $\HCaff$ was introduced earlier by Nazarov
\cite{Naz} (called affine Sergeev algebra)  to study
representations of $\C S_n^-$.

It is known from the works of Sergeev, J\'ozefiak and Stembridge
that the representation theory of $\C S_n^-$ is essentially
equivalent to that of $\Cl_n \rtimes \C S_n$. This phenomenon has
subsequently been clarified by the construction of a superalgebra
isomorphism between $\Cl_n \rtimes \C S_n$ and $\Cl_n \otimes \C
S_n^-$, due to Sergeev \cite{Ser} and Yamaguchi \cite{Yam}
independently. (We will say that $\Cl_n \rtimes \C S_n$ and $\C
S_n^-$ are {\em Morita super-equivalent}; for a justification of
the terminology, cf. \cite[Lemma~9.9]{BK} or \cite[13.2]{Kle}, or
our Section~\ref{sec:morita}). Such a super-equivalence has been
extended by the author \cite{Wa} to one between the  degenerate
spin Hecke algebra introduced in {\em loc. cit.} and Nazarov's
degenerate affine Hecke-Clifford algebra.

%
%
\subsection{Properties of the spin (affine) Hecke algebra}

We establish a Morita super-equivalence between $\HCn$ and $\Hn$
(respectively, between $\HCaff$ and $\Haff$) by constructing
explicitly a $q$-deformed version of the above Morita
super-equivalences \cite{Ser, Yam, Wa} in both finite and affine
setups:
$$\Phi: \HCn \stackrel{\simeq}{\longrightarrow} \Cl_n \otimes
\Hn, \qquad
\Phi: \HCaff \stackrel{\simeq}{\longrightarrow} \Cl_n \otimes
\Haff.$$
Our key observation on the existence of natural subalgebras of
$\HCn$ and $\HCaff$ which super-commute with $\Cl_n$ paves the way
for the presentations of $\Hn$ and $\Haff$.

A fundamental construction in the classical theory of the spin
symmetric group is an algebra homomorphism from $\C S_n^-$ to
$\Cl_{n-1}$ which gives rise to the basic spin $\C
S_n^-$-supermodule \cite{Sch, Joz}. We obtain a natural
$q$-deformation of this construction in which the spin Hecke
algebra $\Hn$ fits nicely.

We construct standard bases for $\Hn$ and for $\Haff$, describe the
center of $\Haff$, and further introduce the intertwiners for
$\Haff$, which have their counterparts in \cite{JN}. We introduce
the cyclotomic spin Hecke algebras and show that they are Morita
super-equivalent to the cyclotomic Hecke-Clifford algebras
introduced in \cite{BK}. We remark that all of the definitions and
constructions in this paper can make sense over a field of
characteristic different from $2$ (which is occasionally assumed to
contain $\sqrt{2}$), and often even over the ring $\Z[\hf]$. It is
possible to develop the representation theory of $\Hn$ and $\Haff$
parallel to the principal results for $\HCn$ and $\HCaff$ in
\cite{BK, JN}. The new perspective of spin (affine) Hecke algebras
can in turn help to clarify the work on the (affine) Hecke-Clifford
algebras.

\subsection{Relation to the (affine) Hecke algebra}
There is a different setup where the spin (affine) Hecke algebra
appears to be relevant. One easily writes down a nonstandard
presentation for the usual Hecke algebra ${\mathcal H}_n$ with new
generators $\tprime_i :=T_i +T_i^{-1}$ instead of the familiar ones
$T_i$. The definition of $\Hn$ and the nonstandard presentation of
${\mathcal H}_n$ are surprisingly compatible and this leads to a
notion of a covering Hecke algebra $\tHn$ which is a $q$-deformation
of $\C \tS_n$.  The quotient of the algebra $\tHn$ by the ideal
$\langle z+1 \rangle$ (respectively, $\langle z-1 \rangle$) is
isomorphic to $\Hn$ (respectively, ${\mathcal H}_n$).

It is remarkable that such a compatibility extends to the spin
affine Hecke algebra $\Haff$ and the usual affine Hecke algebra
$\affH$ of type $GL$, where we have to adopt a nonstandard
presentation of $\affH$ via the generators $\tprime_i$ and $\hf(X_i
\pm X_i^{-1})$ instead of the Bernstein-Lusztig presentation via the
generators $T_i$ and $X_i$. This leads to the definition of the
covering affine Hecke algebra $\tHaff$, whose quotient by the ideal
$\langle z+1 \rangle$ (respectively, $\langle z-1 \rangle$) is
isomorphic to $\Haff$ (respectively, $\affH$).

\subsection{The organization and acknowledgment}

The paper is organized as follows. In Section~\ref{sec:finite}, we
recall the Hecke-Clifford algebra and introduce the spin and
covering Hecke algebras of finite type. In
Section~\ref{sec:structurefinite}, we establish the Morita
super-equivalence between $\HCn$ and $\Hn$. We provide standard
bases for $\Hn$, $\tHn$ and their even subalgebras (which are
$q$-deformations of a double cover of the alternating group and its
spin quotient), and also construct the basic spin $\Hn$-supermodule.
In Section~\ref{sec:affine}, we present the affine Hecke algebra
counterpart of Section~\ref{sec:finite}. In
Section~\ref{sec:structureaffine}, we establish the  Morita
super-equivalence between $\HCaff$ and $\Haff$. We describe the
intertwiners, a standard basis, and the center for $\Haff$. In
Section~\ref{sec:cyclot}, we introduce the cyclotomic spin Hecke
algebras and the Jucys-Murphy elements for $\Hn$. We explain the
degeneration of $\Haff$ and its cyclotomic version. It is our view
that a general notion of spin Hecke algebras exists naturally beyond
the setup in this paper. We will return to this elsewhere.

This research is partially supported by NSF and NSA grants. I
gratefully acknowledge the UVa Sesquicentennial Associateship which
allowed me to spend the Spring semester of 2006 at MSRI, Berkeley. I
thank MSRI for its excellent working atmosphere, where this work was
initiated.

\section{Spin and covering Hecke algebras (of finite type)}
\label{sec:finite}

\subsection{The Hecke-Clifford algebra} Let $q$ be a formal
parameter.
\begin{definition} \cite{Ol}
The {\em Hecke-Clifford algebra} $\HCn$ is the $\C(q)$-algebra
generated by $T_i $ $(1\le i \le n-1)$ and $c_i $ $(1\le i \le
n)$, subject to the following relations:
\begin{eqnarray}
(T_i -q)(T_i +q^{-1}) =0 && \label{hecke}\\
 T_i T_{i+1} T_i = T_{i+1} T_i T_{i+1},
 && T_i T_j = T_j T_i \; (|i-j|>1) \label{braid}\\
T_i c_i = c_{i+1} T_i, && T_i c_j = c_j T_i \;\; (j\neq i,i+1) \label{tici}\\
c_i^2 =1, && c_i c_j = -c_j c_i \; (i \neq j). \label{clifford}
\end{eqnarray}
\end{definition}
The algebra $\HCn$ was introduced by Olshanski \cite{Ol}. It is
naturally a super (i.e. $\Z_2$-graded) algebra with $c_i $ $(1 \le i \le n)$ being odd
and $T_i $ $(1 \le i \le n-1)$ being even. The subalgebra
generated by $T_i$ subject to the relations (\ref{hecke})--(\ref{braid})
is the usual Hecke algebra $\mathcal H_n$ associated
to the symmetric group $S_n$. We define $T_\sigma := T_{i_1}\cdots
T_{i_r}$ as usual for any reduced expression $\sigma
=s_{i_1}\cdots s_{i_r} \in S_n$. The $\C$-algebra generated by
$c_1, \ldots, c_n$ is a Clifford (super)algebra and will be
denoted by $\Cl_n$. It is known that $T_\sigma c_1^{\epsilon_1}
\cdots c_n^{\epsilon_n}$, where $\sigma \in S_n$ and $\epsilon_1,
\cdots, \epsilon_n \in \{0, 1\}$, is a basis for $\HCn$.
Here are some useful identities derived from
(\ref{hecke})--(\ref{clifford}):
\begin{eqnarray}
T_i c_{i+1} &=&c_i T_i - \ep (c_i -c_{i+1}) \label{ticiplus1}\\
T_i (c_i -c_{i+1}) T_i &=& c_{i+1} -c_i   \label{skew} \\
(c_i -c_j)(c_j-c_k)(c_i-c_j) &=& 2(c_k -c_i) \quad \text{for
distinct } i,j,k.  \label{ijk}
\end{eqnarray}

\subsection{The spin Hecke algebra}
Recall
$\ep =q -q^{-1}.$ We now introduce the first new concept of the
paper.
\begin{definition}
The {\em spin Hecke algebra} $\Hn$ is a $\C(q)$-algebra generated
by $\rr_i, 1 \le i \le n-1$, subject to the following relations:
\begin{eqnarray}
\rr_i^2 &=& -\ep^2-2 \equiv -(q^2 +q^{-2}) \label{eq:rri2}\\
\rr_i \rr_j &=& -\rr_j \rr_i  \quad (|i-j|>1)  \label{eq:rrij}\\
\rr_i \rr_{i+1} \rr_i - \rr_{i+1} \rr_i \rr_{i+1} &=& \ep^2
(\rr_{i+1} -\rr_{i}).  \label{braidspin}
\end{eqnarray}
\end{definition}
The algebra $\Hn$ is naturally a superalgebra by requiring each
$\rr_i $ to be odd, since the defining relations for $\Hn$ are
$\Z_2$-homogeneous with respect to such a grading.

\subsection{(Anti-)involutions of $\Hn$}
\label{subsec:invol}

There are several involutions (i.e. algebra automorphisms of order
$2$) of the algebra $\Hn$. Define
\begin{eqnarray*}
\sigma: & \rr_i \mapsto \rr_{n-i}, & q  \mapsto q, \\
s:& \rr_i \mapsto -\rr_i,  &q  \mapsto q, \\
-:&   \rr_i \mapsto \rr_i,   & q \mapsto q^{-1},
\end{eqnarray*}
where $1\le i \le n-1$.
By inspection of the defining relations for $\Hn$, $\sigma, s$ and
$-$ can be extended to homomorphisms of $\Hn$ and they are indeed
involutions of $\Hn$ (regarded as an algebra over $\C$),
Furthermore, $\sigma, s, -$ commute with each other, and their
products give rise to several more involutions.

We also define an anti-involution $\tau$ of $\Hn$ by letting
$\tau (\rr_i) =-\rr_i$ for each $i.$ One obtains more
anti-involutions by composing $\tau$ with the involutions above
(which commute with $\tau$).
\subsection{A nonstandard presentation of Hecke algebra}

Denote
$$\tprime_i: =T_i +T_i^{-1} \equiv 2T_i -\ep.$$

\begin{proposition}
The algebra ${\mathcal H}_n$ has a presentation with generators $\tprime_i$ $(1 \le i \le n-1)$
subject to the following relations:
\begin{eqnarray}
\tprime_i^2 &=& q^2 +q^{-2} +2  \label{eq:tprime}\\
\tprime_i \tprime_j &=& \tprime_j \tprime_i \quad (|i-j|>1) \label{eq:titjprime}\\
\tprime_i \tprime_{i+1} \tprime_i - \tprime_{i+1} \tprime_i
\tprime_{i+1} &=& \ep^2 (\tprime_{i+1}
-\tprime_i).\label{eq:tiplus1prime}
\end{eqnarray}
\end{proposition}
\begin{proof}
Follows by a direct computation.
%
\end{proof}


%
%
\subsection{The covering Hecke algebra}
\begin{definition}
The {\em covering Hecke algebra} $\tHn$ is a $\C(q)$-superalgebra
generated by the even generator $z$ and the odd generators
$\tilde{T}_i$ $(1 \le i \le n-1)$, subject to the following
relations:
\begin{eqnarray}
z^2=1, & & z \text{ is central} \nonumber \\
\tilde{T}_i^2 &=& z(q^2 +q^{-2} +1) +1 \label{eq:zt}\\
\tilde{T}_i \tilde{T}_j &=& z\tilde{T}_j \tilde{T}_i \quad (|i-j|>1) \label{eq:ztitj}\\
\tilde{T}_i \tilde{T}_{i+1} \tilde{T}_i - \tilde{T}_{i+1}
\tilde{T}_i \tilde{T}_{i+1} &=& \ep^2 (\tilde{T}_{i+1}
-\tilde{T}_{i}). \label{eq:ztitiplus1}
\end{eqnarray}
\end{definition}
We shall denote by $\langle a, b, \cdots \rangle$ a two-sided ideal
generated by $a, b, \cdots.$ The quotient of the covering Hecke
algebra $\tHn$ by the ideal $\langle z-1\rangle$ is isomorphic to
the usual Hecke algebra $\mathcal H_n$ with nonstandard presentation
(where the canonical image of $\tilde{T}_i$ matches $\tprime_i$) and
the quotient by $\langle z+1 \rangle$ is isomorphic to the spin
Hecke algebra $\Hn$ (where the canonical image of $\tilde{T}_i$
matches $\rr_i$).

\section{Algebraic structures of the spin Hecke algebra}
\label{sec:structurefinite}
\subsection{A Morita super-equivalence}
\label{sec:morita}

Note that the multiplication in a tensor product $\Cl \otimes \Bn$
of two superalgebras $\Cl$ and $\Bn$ has a suitable sign
convention:
$$(c' \otimes b')(c \otimes b) =(-1)^{|b'||c|} (c'c
\otimes b'b).$$
We shall write a typical element in $\Cl \otimes
\Bn$ as $cb$ rather than $c \otimes b$, and use short-hand
notations $c =c \otimes 1, b =1\otimes b$.


\begin{theorem}  \label{th:isomfiniteq}
There exists a superalgebra isomorphism
\begin{eqnarray*}
\Phi: \HCn \stackrel{\simeq}{\longrightarrow} \Cl_n \bigotimes \Hn
\end{eqnarray*}
which extends the identity map on $\Cl_n$ and sends
\begin{equation} \label{phiT}
T_i \mapsto T_i^\Phi
 := -\hf \rr_i (c_i -c_{i+1}) +\frac{\ep}2 (1 -c_ic_{i+1}), \quad
1\le i \le n-1.
\end{equation}
Its inverse map $\Psi$ extends the identity map on $\Cl_n$ and
sends
\begin{equation} \label{psirr}
 \rr_i \mapsto \rr_i^\Psi := (c_i -c_{i+1}) T_i +\ep c_{i+1},
 \quad  1\le i \le n-1.
\end{equation}
\end{theorem}

\begin{remark}
The isomorphism in Theorem~\ref{th:isomfiniteq}
in the $q \mapsto 1$ limit reduces to the
superalgebra isomorphism $\Cl_n \rtimes \C S_n \cong \Cl_n \otimes \C S_n^-$
found in \cite{Ser, Yam}.
\end{remark}

By (\ref{hecke})--(\ref{clifford}), we have the following
equivalent expressions for $ \rr_i^\Psi$:
\begin{eqnarray}
 \rr_i^\Psi &=&  -T_i(c_i -c_{i+1}) +\ep c_{i} \nonumber
\\
&=& c_i T_i -c_{i+1} T_i^{-1} = T_i c_{i+1}  -T_i^{-1} c_i.
\label{equivexp}
\end{eqnarray}

Thanks to the above isomorphisms $\Phi, \Psi$, we can define exact
functors
\begin{eqnarray*}
\mathfrak F: \Hn\text{-smod} \rightarrow \HCn\text{-smod},
 &&
\mathfrak F: = \Phi^* (U_n  \otimes  ?), \\
\mathfrak G: \HCn\text{-smod} \rightarrow \Hn\text{-smod},
 && \mathfrak G := \text{Hom}_{\Cl_n} (U_n, \Psi^*(?)),
\end{eqnarray*}
where $U_n$ denotes the basic spin $\Cl_n$-supermodule and
$\Hn$-smod (respectively, $\HCn$-smod) denotes the category of
finite-dimensional supermodules of $\Hn$ (respectively, of $\HCn$).
For $n$ even, $\mathfrak F$ and $\mathfrak G$ establish the
equivalence of categories, and indeed $\Hn$ and $\HCn$ are Morita
equivalent in the usual sense since $\Cl_n$ is a simple algebra. For
$n$ odd, $\Cl_n$ is a simple {\em super}algebra of type $Q$,
$\mathfrak F$ and $\mathfrak G$ establish an {\em almost} Morita
equivalence of categories which involves some $\Z_2$-parity change
functor (see \cite[Lemma~9.9]{BK} or \cite[Proposition~13.2.2]{Kle}
for a precise statement in a similar setup).

Let us call two superalgebras $\An$ and $\Bn$ {\em Morita
super-equivalent} if there is a superalgebra isomorphism $\An
\cong \Cl_n \otimes \Bn$ or $\Bn \cong \Cl_n \otimes \An$ for some
Clifford algebra $\Cl_n$.  In particular, $\HCn$ and $\Cl_n
\otimes \Hn$ are Morita super-equivalent. (This restrictive
definition of the Morita super-equivalence is all we need in this
paper, though it is possible and potentially useful in other
contexts to give a more general definition which incorporates the
usual Morita equivalence.)

\subsection{Proof of the isomorphism Theorem~\ref{th:isomfiniteq}}

We start with several lemmas.
\begin{lemma} \label{lem:commute}
The $ \rr_i^\Psi$ super-commute with $\Cl_n$ for every $1\le i \le
n-1$.
\end{lemma}
\begin{proof}
Clearly, $\rr_i^\Psi c_j = -c_j \rr_i^\Psi$ for $j \neq i,i+1$. By
(\ref{psirr}), (\ref{tici}) and (\ref{clifford}),
$$\rr_i^\Psi c_i
 = (c_i -c_{i+1})c_{i+1} T_i +\ep  c_{i+1} c_i
 = -c_i(c_i -c_{i+1}) T_i -\ep c_ic_{i+1}
 = -c_i \rr_i^\Psi.$$
We leave to the reader the similar verification that $\rr_i^\Psi
c_{i+1} =-c_{i+1}  \rr_i^\Psi$.
\end{proof}

\begin{lemma}  \label{lem:homom}
The $ \rr_i^\Psi$, $1\le i \le n-1$, satisfy the relations
(\ref{eq:rri2}) and (\ref{eq:rrij}).
\end{lemma}
\begin{proof}
Clearly, $\rr_i^\Psi \rr_j^\Psi  = -\rr_j^\Psi\rr_i^\Psi$ for
$|i-j|>1.$
By (\ref{hecke})--(\ref{clifford}) and (\ref{equivexp}), we
calculate that
\begin{eqnarray*}
(\rr_i^\Psi)^2 &=&
(c_i T_i -c_{i+1} T_i^{-1})(T_i c_{i+1}  -T_i^{-1} c_i)  \\
 &=& c_i (\ep T_i +1) c_{i+1} +c_{i+1} ( 1+\ep^2 -\ep T_i) c_i -2 \\
 &=& c_i (\ep T_i^{-1} +1 +\ep^2) c_{i+1} +c_{i+1} ( 1+\ep^2 -\ep T_i) c_i -2 \\
 &=& c_i \ep T_i^{-1} c_{i+1} -c_{i+1} \ep T_i c_i -2
 =-\ep^2-2.
\end{eqnarray*}
This verifies (\ref{eq:rri2}) and (\ref{eq:rrij}) for
$\rr_i^\Psi$.
\end{proof}

\begin{lemma} \label{lem:braidmod}
The relation (\ref{braidspin}) holds for $ \rr_i^\Psi$, that is,
for every $1\le i \le n-2$,
 $$\left( \rr_i^\Psi\rr_{i+1}^\Psi  +\ep^2 \right)\rr_i^\Psi
 = \rr_{i+1}^\Psi
 \left( \rr_i^\Psi\rr_{i+1}^\Psi  +\ep^2 \right).
 \label{eq:rriiplus1ipsi}$$
\end{lemma}
\begin{proof}
Formulas (\ref{hecke}) through (\ref{ijk}) are frequently used in
this proof. We have
\begin{eqnarray*}
\rr_i^\Psi \rr_{i+1}^\Psi
 &=& ( (c_i-c_{i+1}) T_i^{-1} +\ep c_i)
      ((c_{i+1} -c_{i+2}) T_{i+1} +\ep c_{i+2}) \\
 &=& (c_i -c_{i+1})(c_i -c_{i+2}) (T_i -\ep) T_{i+1}
  +\ep c_i (c_{i+1} -c_{i+2}) T_{i+1}\\
  && +\ep (c_i -c_{i+1}) c_{i+2} (T_i -\ep) +\ep^2 c_i c_{i+2} \\
 &=& (c_i -c_{i+1})(c_i -c_{i+2}) T_i T_{i+1}
  +\ep c_{i+2} (c_{i+1} -c_{i+2}) T_{i+1} \\
  && -\ep  c_{i+2} (c_i -c_{i+1}) T_i
   +\ep^2 c_{i+1} c_{i+2}.
\end{eqnarray*}

Recalling $\rr_i^\Psi = (c_i -c_{i+1}) T_i +\ep c_{i+1}$ from
(\ref{psirr}), we have that
\begin{eqnarray*}
&&\left( \rr_i^\Psi\rr_{i+1}^\Psi  +\ep^2 \right)\rr_i^\Psi \\
 &=& (c_i -c_{i+1})(c_i -c_{i+2})(c_{i+1} -c_{i+2}) T_iT_{i+1}T_i  \\
  &&
  +\ep c_{i+2} (c_{i+1} -c_{i+2})(c_i -c_{i+2}) T_{i+1}T_i \\
  && -\ep  c_{i+2} (c_i -c_{i+1})(c_{i+1} -c_{i})
   +\ep^2 (c_i-c_{i+1} +c_{i+2}+c_ic_{i+1}c_{i+2}) T_i  \\
 && +\ep(c_i -c_{i+1})(c_i -c_{i+2}) c_{i+2} T_i T_{i+1}
  +\ep^2 c_{i+2} (c_{i+1} -c_{i+2}) c_{i+2} T_{i+1} \\
  && -\ep^2  c_{i+2} (c_i -c_{i+1})c_i T_i
   +\ep^3 c_{i+2} (c_i -c_{i+1}) (c_i -c_{i+1})
   +\ep^3 (c_{i+1} -c_{i+2}) \\
 &=& 2(c_i -c_{i+2}) T_iT_{i+1}T_i
  +\ep(c_{i+1} +c_{i+2} -c_i -c_ic_{i+1}c_{i+2}) T_{i+1}T_i \\
  &&
  +\ep(c_{i+1} +c_{i+2} -c_i +c_ic_{i+1}c_{i+2}) T_iT_{i+1}
   -\ep^2 (c_{i+1} +c_{i+2}) T_{i+1}\\
  &&  +\ep^2 (c_i -c_{i+1}) T_i
  +2\ep c_{i+2} +\ep^3 (c_{i+1} +c_{i+2}).
\end{eqnarray*}

On the other hand, recalling $\rr_{i+1}^\Psi
=T_{i+1}(c_{i+2}-c_{i+1}) +\ep c_{i+1}$, we have
\begin{eqnarray*}
&& \rr_{i+1}^\Psi \left( \rr_i^\Psi\rr_{i+1}^\Psi  +\ep^2 \right) \\
 &=& 2(c_i -c_{i+2}) T_{i+1}T_iT_{i+1}
  +\ep T_{i+1}(c_{i+2}-c_{i+1})c_{i+2} (c_{i+1} -c_{i+2}) T_{i+1} \\
  && -\ep T_{i+1}(c_{i+2}-c_{i+1}) c_{i+2} (c_i -c_{i+1}) T_i
   -2\ep^2 c_{i+2}T_{i+1} \\
  && +\ep c_{i+1}(c_i -c_{i+1})(c_i -c_{i+2}) T_i T_{i+1}
  +\ep^2 c_{i+1} c_{i+2} (c_{i+1} -c_{i+2}) T_{i+1} \\
  && -\ep^2 c_{i+1}  c_{i+2} (c_i -c_{i+1}) T_i
   +\ep^3 (c_{i+1} +c_{i+2}) \\
 &=& 2(c_i -c_{i+2}) T_{i+1}T_iT_{i+1}
  + 2\ep c_{i+2} (\ep T_{i+1}+1) \\
  && -\ep (c_{i+1}-c_{i+2}) c_{i+1} (c_i -c_{i+2}) T_{i+1}T_i
   +\ep^2 (c_{i+1}-c_{i+2}) c_{i+1} (c_i -c_{i+1}) T_i \\
  && -2\ep^2 c_{i+2}T_{i+1}
  +\ep c_{i+1}(c_i -c_{i+1})(c_i -c_{i+2}) T_i T_{i+1}
  +\ep^2 c_{i+1} c_{i+2} (c_{i+1} -c_{i+2}) T_{i+1} \\
  && -\ep^2 c_{i+1}  c_{i+2} (c_i -c_{i+1}) T_i
   +\ep^3 (c_{i+1} +c_{i+2}),
\end{eqnarray*}
which can be shown by a simple rewriting to coincide with the
right-hand side of the previous formula for
$(\rr_i^\Psi\rr_{i+1}^\Psi +\ep^2 )\rr_i^\Psi$.
\end{proof}

\begin{lemma}  \label{lem:inverse}
For every $1\le i \le n-1$, we have $\Psi (T_i^\Phi) =T_i$, and
$\Phi (\rr_i^\Psi) =\rr_i$.
\end{lemma}
\begin{proof}
By (\ref{equivexp}), we have
\begin{eqnarray*}
\Psi (T_i^\Phi) &=&  \Psi \left(\hf \rr_i (c_{i+1} -c_i)
+\frac{\ep}2 (1 -c_ic_{i+1}) \right) \\
&=& \hf (-T_i(c_i -c_{i+1}) +\ep c_{i})(c_{i+1} -c_i)
 + \hf {\ep} (1 -c_ic_{i+1}) =T_i, \\
\Phi(\rr_i^\Psi) &=& \Phi (T_i(c_{i+1} -c_i) +\ep c_{i}) \\
&=& -\left(\hf \rr_i (c_{i+1} -c_i) +\frac{\ep}2 (1 -c_ic_{i+1})
\right) (c_i -c_{i+1}) +\ep c_{i}  =\rr_i.
\end{eqnarray*}
\end{proof}

\begin{proof}[Proof of Theorem~\ref{th:isomfiniteq}]
By Lemmas~\ref{lem:commute}, \ref{lem:homom}, \ref{lem:braidmod},
$\Psi$ is a (super) algebra homomorphism. By
Lemma~\ref{lem:inverse} and $\Psi (c_i) =c_i$, $\Psi$ is
surjective.

Denote by ${\mathcal H}_{n,\Psi}^{-}$ the subalgebra of $\HCn$
generated by $\rr_i^\Psi, 1\le i \le n-1$. By
Lemma~\ref{lem:commute}, we have
$\HCn \supseteq \Cl_n \otimes {\mathcal H}_{n,\Psi}^{-}.$ By
Lemma~\ref{lem:inverse}, we have
$$T_i = \hf  (c_i -c_{i+1}) \rr_i^\Psi +\frac{\ep}2
(1 -c_ic_{i+1}) \in \Cl_n \otimes {\mathcal H}_{n,\Psi}^{-},$$
and thus all generators $T_i, c_i$ of $\HCn$ lie in $\Cl_n \otimes
{\mathcal H}_{n,\Psi}^{-}$. Therefore, $\HCn = \Cl_n \otimes
{\mathcal H}_{n,\Psi}^{-}$ and $\dim {\mathcal H}_{n,\Psi}^{-}
=n!$. By Proposition~\ref{prePBW} below (whose proof is elementary
and in particular does not use this Theorem), we have $\dim \Hn
\le n!$. Thus for dimension reason the surjective homomorphism
$\Psi |_{\Hn}: \Hn \rightarrow {\mathcal H}_{n,\Psi}^{-}$ is
indeed an isomorphism and $\dim \Hn =n!$.
Since both $\HCn$ and $\Cl_n \otimes \Hn$ have dimensions equal to
$2^n n!$, the surjective homomorphism $\Psi$ is an algebra
isomorphism.

By Lemma~\ref{lem:inverse}, $\Psi$ and $\Phi$ are inverse
isomorphisms.
\end{proof}

\begin{remark}
A somewhat different argument for Theorem~\ref{th:isomfiniteq} goes
as follows. We can verify directly that $\Phi$ is an algebra
homomorphism in a way similar to $\Psi$, which involves a tedious
verification of the braid relations for $T_i^\Phi$. Then
Lemma~\ref{lem:inverse} implies that $\Phi$ and $\Psi$ are inverse
isomorphisms. This argument does not use Proposition~\ref{prePBW}
below.
\end{remark}

\subsection{Bases for $\Hn$ and $\tHn$}

Introduce the following monomials in $\Hn$:
\begin{eqnarray} \label{monomial}
\rr_{i,a} :=\rr_i\rr_{i-1} \cdots \rr_{i-a+1}, \quad 0\le a \le i,
\; 1\le i \leq n-1,
\end{eqnarray}
where it is understood that $\rr_{i,0} \equiv 1$ for all $i$.
%
We refer to a product $\rr_{i_1}\rr_{i_2} \cdots \rr_{i_s}$ of
generators in $\Hn$ as a {\em monomial} in $\Hn$, and call a
monomial {\em standard} if it is of the form
$\rr_{1,a_1}\rr_{2,a_2}\cdots \rr_{n-1,a_{n-1}}$, where $0 \le
a_i\le i$ for each $1 \le i \le n-1$.

\begin{proposition}  \label{prePBW}
The standard monomials $\rr_{1,a_1}\rr_{2,a_2}\cdots
\rr_{n-1,a_{n-1}}$, where $0 \le a_i\le i$ and $1 \le i \le n-1$,
linearly span $\Hn$. In particular, $\dim \Hn \le n!$.
\end{proposition}

\begin{proof}
Since the number of standard monomials in $\Hn$ is $n!$, it
suffices to prove the first statement on linear span.

{\em Claim 1.} Any monomial of $\Hn$ is spanned by the monomials
in which $\rr_{n-1}$ appears at most once.

We prove Claim~1 by induction on $n$. The claim trivially holds for
$n=1$. For any monomial in which $\rr_{n-1}$ appears more than
twice, we can apply the argument below to a portion of the monomial
which contains exactly two $\rr_{n-1}$'s to reduce the number of
$\rr_{n-1}$'s. So, let us assume that a given monomial is of the
form $\rr_{n-1} \cdot \rr_{i_1} \rr_{i_2} \cdots \rr_{i_s}\cdot
\rr_{n-1}$, where $1\le i_1, \cdots, i_s \le n-2$ (for some $s$). By
applying (\ref{eq:rrij}) to move the $\rr_i$'s outbound the two
$\rr_{n-1}$'s whenever possible, we are reduced to the case $s=0$ or
${i_1} ={i_s} =n-2$, where $s \ge 1$. The reduction of the number of
$\rr_{n-1}$'s in the case $s=0$ is done by (\ref{eq:rri2}), while in
the case $s=1$ is by (\ref{braidspin}) (note that here we got a
linear combination of monomials because (\ref{braidspin}) is not the
usual braid relation). In the case $s \ge 2$, we reduce to the
previous cases by applying the claim for $n-1$, which is the
induction step.

{\em Claim 2.} Any monomial of $\Hn$ in which $\rr_{n-1}$ appears
exactly once can be written as a linear combination of monomials
of the form $\rr_{i_1}\rr_{i_2} \cdots \rr_{i_s}\cdot
\rr_{n-1,a_{n-1}}$, where $1\le a_{n-1} \le n-1$ and $1\le i_1,
\ldots, i_s \le n-2$ for some $s$.

We again argue by induction on $n$. The claim is trivial for
$n=1,2$. By permuting $\rr_{n-1}$ via (\ref{eq:rrij}) to the right
as much as possible, we rewrite the monomial (up to a sign) such
that $\rr_{n-1}$ appears at the very end or it is followed by
$\rr_{n-2}$. We continue with the second possibility, otherwise we
are done. Since the part of the monomial starting from $\rr_{n-2}$
to the right lies in $H^-_{n-1}$, we may apply Claim~1 to reduce to
the case when $\rr_{n-2}$ appears to the right of $\rr_{n-1}$
exactly once. Now the induction step for $n-1$ of Claim~2 is
applicable to complete the proof of Claim~2.

We now proceed by induction on $n$. The proposition holds
trivially for $n=1$. If a monomial in $\Hn$ does not contain
$\rr_{n-1}$ and thus is a monomial in ${\mathcal H}_{n-1}^-$, then
it is a linear combination of the standard monomials as the
induction step applies. Otherwise, the proposition follows by the
induction step, Claims~1 and 2.
\end{proof}

\begin{theorem} \label{th:PBWspin}
The standard monomials $\rr_{1,a_1}\rr_{2,a_2}\cdots
\rr_{n-1,a_{n-1}}$, where $0 \le a_i\le i$ and $1 \le i \le n-1,$
form a basis for $\Hn$. Also, $\dim \Hn =n!$.
\end{theorem}
\begin{proof}
The statement that $\dim \Hn =n!$ is a consequence of (the proof
of) Theorem~\ref{th:isomfiniteq}. The number of standard monomials
in $\Hn$ is $n!$, and thus the theorem follows from
Proposition~\ref{prePBW}.
\end{proof}

We define the monomial $\tT_{i,a} \in \tHn$ (respectively,
$\tprime_{i,a} \in \mathcal H_n$), with $\tT$ (respectively,
$\tprime$) replacing $\rr$ in the definition (\ref{monomial}) of
$\rr_{i,a}$.

\begin{proposition} \label{PBWcovering}
The elements $\tT_{1,a_1}\tT_{2,a_2}\cdots \tT_{n-1,a_{n-1}},$
$z\tT_{1,a_1}\tT_{2,a_2}\cdots \tT_{n-1,a_{n-1}}$, where $0 \le
a_i\le i$ and $1 \le i \le n-1,$ form a basis for $\tHn$. Also,
$\dim \tHn =2n!$.
\end{proposition}

\begin{proof}
The same argument for Prop.~\ref{prePBW} shows that the elements
in the proposition form a spanning set for $\tHn$. It remains to
prove the linear independence.

By definitions, $\tHn / \langle z-1\rangle \cong \mathcal H_n$ and
$\tHn / \langle z+1 \rangle \cong \Hn$. Denote the corresponding
canonical maps by $p_+: \tHn \rightarrow \mathcal H_n$ and $p_-:
\tHn \rightarrow \Hn$. Clearly, $p_+(\tT_i) =\tprime_i$,
$p_-(\tT_i) =\rr_i$, and $p_{\pm} (z) =\pm 1$.
We shall use the short-hand notation $\tT_{\bf a}
=\tT_{1,a_1}\tT_{2,a_2}\cdots \tT_{n-1,a_{n-1}}$. Assume there is
a relation $(\star)$
$\sum_{\bf a} (\alpha_{\bf a} \tT_{\bf a} +  \beta_{\bf a} z
\tT_{\bf a}) =0$ for some constants $\alpha_{\bf a}, \beta_{\bf
a}.$ By applying the canonical map $p_-$ to $(\star)$ and
Theorem~\ref{th:PBWspin}, we conclude that $\alpha_{\bf a} -
\beta_{\bf a}=0$ for each $\bf a$. On the other hand, it is (a
variant of) a classical fact that
$\tprime_{1,a_1}\tprime_{2,a_2}\cdots \tprime_{n-1,a_{n-1}}$,
where $0 \le a_i\le i$ and $1 \le i \le n-1$, form a linear basis
for $\Hn$. By applying $p_+$ to $(\star)$, we conclude that
$\alpha_{\bf a} +\beta_{\bf a}=0$ for each $\bf a$. So
$\alpha_{\bf a} = \beta_{\bf a} \equiv 0$.
\end{proof}

\begin{remark}
By (\ref{eq:rri2})--(\ref{braidspin}) and
Theorem~\ref{th:PBWspin}, $\Hn$ reduces to $\C S_n^-$ as $q$ goes
to $1$. By a standard deformation argument, the algebra $\Hn$ is
semisimple. Similarly, $\tHn$ is a flat deformation of $\C \tS_n$
by Proposition~\ref{PBWcovering}, which justifies the terminology
of ``covering Hecke algebra" for $\tHn$.
\end{remark}
\subsection{Spin Hecke algebra for the alternating group}

Note that the number of generators $\rr_i$'s appearing in the
monomial $\rr_{1,a_1}\rr_{2,a_2}\cdots \rr_{n-1,a_{n-1}}$ is $a_1
+\ldots +a_{n-1}$ (called the {\em length} of the monomial).
Denote by $\mathcal H^-_{n,\bar{0}}$ (respectively,
$H^\thicksim_{n,\bar{0}}$) the even subalgebra of the superalgebra
$\Hn$ (respectively, $\tHn$).

\begin{proposition}
Let $n \ge 2$.

\begin{enumerate}
\item The standard monomials of even length, that is,
$\rr_{1,a_1}\rr_{2,a_2}\cdots \rr_{n-1,a_{n-1}}$, where $0 \le
a_i\le i,$ $1 \le i \le n-1$ such that $a_1 +\ldots +a_{n-1}$ is
even, form a basis for the algebra $\mathcal H^-_{n,\bar{0}}$. In
particular, $\dim \mathcal H^-_{n,\bar{0}} =\hf n!$.
%
 \item
The elements $\tT_{1,a_1}\tT_{2,a_2}\cdots \tT_{n-1,a_{n-1}},
z\tT_{1,a_1}\tT_{2,a_2}\cdots \tT_{n-1,a_{n-1}}$, where $0 \le
a_i\le i,$ $1 \le i \le n-1$ such that $a_1 +\ldots +a_{n-1}$ is
even, form a basis for the algebra $\mathcal
H^\thicksim_{n,\bar{0}}$. In particular, $\dim \mathcal
H^\thicksim_{n,\bar{0}} =n!$.
\end{enumerate}
\end{proposition}

\begin{proof}
Follows from Theorem~\ref{th:PBWspin} and
Proposition~\ref{PBWcovering}.
\end{proof}

Recall that the alternating group $A_n$ is a subgroup of $S_n$ of
index $2$. The short exact sequence (\ref{centralext}) gives rise
to a subgroup $A^\thicksim_n$ of $\tS_n$ of index $2$ which is a
double cover of $A_n$. It follows that $\mathcal
H^\thicksim_{n,\bar{0}}$ is a $q$-deformation of the algebra $\C
A^\thicksim_n$, while $\mathcal H^-_{n,\bar{0}}$ is a
$q$-deformation of the algebra $\C A^\thicksim_n/\langle
z+1\rangle$.

\begin{definition}
The algebra $\mathcal H^-_{n,\bar{0}}$ is called the {\em spin
Hecke algebra for the alternating group $A_n$}. The algebra
$\mathcal H^\thicksim_{n,\bar{0}}$ is called the {\em covering
Hecke algebra for the alternating group $A_n$}.
\end{definition}

We leave it to the reader to write down a presentation for the
algebra $\mathcal H^-_{n,\bar{0}}$ using the generators
$\rr_1\rr_{i+1}$ $(1\le i \le n-2)$ and a similar presentation for
$\mathcal H^\thicksim_{n,\bar{0}}$.

\subsection{The basic spin supermodule}

\begin{theorem}
There exists a homomorphism of superalgebras
$$\pi_q: \Hn \rightarrow \Cl_n \otimes \C (q)$$
 which sends
$$ \rr_i \mapsto \sqrt{-1} (qc_i -q^{-1} c_{i+1}), \quad 1\le i \le n-1.$$
The image is isomorphic to $\Cl_{n-1} \otimes \C (q)$.
\end{theorem}

\begin{proof}
We need to check the relations (\ref{eq:rri2})--(\ref{braidspin})
with  $\g_i :=\sqrt{-1} (qc_i -q^{-1} c_{i+1})$ replacing  $\rr_i$
therein. Clearly we have
$$\g_i^2 = -2 -\ep^2, \quad \g_i \g_j = -\g_j \g_i \;\;
(|i-j|>1).$$

By a straightforward computation, we have
\begin{eqnarray*}
\g_i \g_{i+1} \g_i
 &=& \sqrt{-1} \left(2q c_i - (q^{-1} -q^3)c_{i+1} -(q^{-3}+q)c_{i+2} \right) \\
 &=& 2 \g_i +(q^2 +q^{-2}) \g_{i+1}, \\
\g_{i+1} \g_i \g_{i+1}
 &=& \sqrt{-1} \left((q^3+q^{-1}) c_i -(q^{-3} -q) c_{i+1} -2q^{-1} c_{i+2}\right) \\
 &=&  (q^2 +q^{-2}) \g_i  +2\g_{i+1}.
\end{eqnarray*}
Thus, $\g_i \g_{i+1} \g_i -\g_{i+1} \g_i \g_{i+1} =\ep^2 (\g_{i+1}
-\g_i).$

Since the image of the linear span of $\rr_i$ $(1\le i \le n-1)$
is by definition a subspace of dimension $n-1$ of the linear span
of $c_i$'s, the image of $\Hn$ under $\pi$ is a Clifford algebra
in $n-1$ generators.
\end{proof}
It is well known that $\Cl_{n-1}$ has a unique simple supermodule,
which is of dimension $2^{[n/2]}$. Here $[n/2]$ denotes the
largest integer no greater than $n/2$. The pullback via $\pi_q$
gives rise to a simple $\Hn$-supermodule of dimension $2^{[n/2]}$,
which we will refer to as the {\em basic spin $\Hn$-supermodule}.
Indeed, this module is a $q$-deformation of the basic spin $\C
S_n^-$-supermodule and the homomorphism $\pi_q$ is the
$q$-deformation of a classical fundamental construction
(\cite{Sch, Joz}).

\section{Spin and covering affine Hecke algebras}
\label{sec:affine}

\subsection{The affine Hecke-Clifford algebra}

\begin{definition} \cite{JN}
The {\em affine Hecke-Clifford algebra} $\HCaff$ is the
$\C(q)$-algebra generated by $T_i $ $(1\le i \le n-1)$ and $c_i,
X_i^{\pm 1} $ $(1\le i \le n)$, subject to the relations
(\ref{hecke})--(\ref{clifford}) of $T_i, c_i$ in $\HCn$ and the
following additional relations:
\begin{eqnarray}
(T_i +\ep c_i c_{i+1}) X_{i} T_i &=& X_{i+1} \label{eq:txt}\\
 T_i X_j &=& X_j T_i \;\;\; (j \neq i, i+1)  \label{eq:tixj} \\
X_iX_j &=& X_jX_i   \label{eq:xx}\\
X_i c_i = c_{i} X_i^{-1}, && X_i c_j = c_j X_i \;\;\;\; (i\neq j).
 \label{eq:xc}
\end{eqnarray}
\end{definition}

The affine Hecke-Clifford algebra $\HCaff$ was introduced by
Jones-Nazarov \cite{JN}. The algebra $\HCaff$ admits a canonical
superalgebra structure with $T_i, X_i$ being even and $c_i$ being
odd. It is known that $X_1^{\alpha_1} \cdots X_n^{\alpha_n}
c_1^{\epsilon_1} \cdots c_n^{\epsilon_n} T_\sigma$, where
$\alpha_1,\ldots,\alpha_n \in \Z$, $\epsilon_1,\ldots,\epsilon_n
\in \{0,1\}$ and $\sigma \in S_n$, form a {\em standard} basis for
$\HCaff$ \cite{JN} (also cf. \cite{BK}). By definition, $\HCaff$
contains $\HCn$ as a subalgebra.

The convention that $c_i^2=-1$ was used in \cite{JN}, and so our
$c_i$ matches with their $\sqrt{-1} c_i$. Our convention that
$c_i^2=1$ is consistent with \cite{BK, Kle}. The different
convention leads to a different sign whenever a quadratic term
$c_ic_j$ appears. The following useful identities follow from the
definition:
\begin{eqnarray*}
 (T_i +\ep c_i c_{i+1})^{-1} &=& T_i +\ep c_i c_{i+1} -\ep  \\
T_i X_i &=& X_{i+1}T_i -\ep (X_{i+1} +c_ic_{i+1} X_i) \\
T_i X_{i+1} &=& X_{i}T_i +\ep (1 -c_ic_{i+1})X_{i+1}.
\end{eqnarray*}

\subsection{The spin affine Hecke algebra}

Now we introduce the main new concept of the paper.
\begin{definition}
The {\em spin affine Hecke algebra}, denoted by $\Haff$, is the
$\C(q)$-algebra generated by $\rr_i $ $(1 \le i \le n-1)$ and
$p_i, q_i $ $(1 \le i \le n)$, subject to the relations
(\ref{eq:rri2})--(\ref{braidspin}) for $R_i$'s in $\Hn$ and the
following additional relations:
\begin{eqnarray}
\pp_i \pp_j = \pp_j\pp_i,
 && \qq_i \qq_j = -\qq_j\qq_i  \;\; (i \neq j) \label{ppqqcommute} \\
\pp_i^2 +\qq_i^2 =1, && \pp_i \qq_j = \qq_j\pp_i  \quad (\forall
i, j)   \label{pq} \\
\rr_i \pp_{j} = \pp_{j}\rr_i, &&
 \rr_i \qq_{j} = -\qq_{j}\rr_i \quad (j\neq i,i+1)
 \label{rpqcommute}
\end{eqnarray}
\begin{eqnarray}
\rr_i \pp_i
 &=& \pp_{i+1}\rr_i +\ep (\qq_i -\qq_{i+1}) \label{eq:rp} \\
\rr_i \qq_i
 &=& -\qq_{i+1}\rr_i -\ep (\pp_i +\pp_{i+1}) \label{eq:rq}
\end{eqnarray}
\end{definition}
The algebra $\Haff$ has a canonical superalgebra structure with
each $\pp_i$ being even and each $\qq_i$, $\rr_i$ being odd.

\begin{proposition}  \label{prop:recur}
Assume only the relation (\ref{eq:rri2}). The three pairs of
relations (\ref{eq:rp})--(\ref{eq:rq}),
(\ref{eq:rpp})--(\ref{eq:rqq}) and
(\ref{eq:recurp})--(\ref{eq:recurq}) are equivalent:
\begin{eqnarray}
\rr_i \pp_{i+1}
 &=&  \pp_{i}\rr_i -\ep (\qq_i -\qq_{i+1})  \label{eq:rpp} \\
\rr_i \qq_{i+1}
 &=& -\qq_{i}\rr_i -\ep (\pp_i +\pp_{i+1});  \label{eq:rqq}
\end{eqnarray}
\begin{eqnarray}
 p_{i+1} &=& -\hf \rr_i p_i \rr_i + \frac{\ep}{2} (q_i \rr_i +\rr_i
 q_i) +\frac{\ep^2}{2} p_i \label{eq:recurp} \\
 q_{i+1} &=&  \hf \rr_i q_i \rr_i + \frac{\ep}{2} (p_i \rr_i +\rr_i
 p_i) -\frac{\ep^2}{2} q_i.  \label{eq:recurq}
\end{eqnarray}
In particular, the algebra $\Haff$ is generated by $p_1, q_1$ and
$\rr_i $ $(1 \le i \le n-1)$.
\end{proposition}

\begin{proof}
The last statement follows readily from
(\ref{eq:recurp})--(\ref{eq:recurq}). Recall from (\ref{eq:rri2})
that $\rr_i^2 =-2-\ep^2$.

(i). (\ref{eq:rp})--(\ref{eq:rq}) $\Rightarrow$
(\ref{eq:recurp})--(\ref{eq:recurq}):
The right multiplication of (\ref{eq:rp}) by $\rr_i$ gives us
\begin{eqnarray} \label{eq:rqr}
\rr_i\qq_i \rr_i =(-2-\ep^2)\pp_{i+1} +\ep\qq_i \rr_i
-\ep\qq_{i+1} \rr_i.
\end{eqnarray}
Rewrite (\ref{eq:rq}) as
$ -\qq_{i+1}\rr_i  = \rr_i \qq_i +\ep (\pp_i +\pp_{i+1}). $
Plugging this equation into (\ref{eq:rqr}) and reorganizing the
terms, we obtain (\ref{eq:recurp}). The proof of (\ref{eq:recurq})
is almost identical.

(ii). (\ref{eq:recurp})--(\ref{eq:recurq})  $\Rightarrow$
(\ref{eq:rpp})--(\ref{eq:rqq}):
Multiplying (\ref{eq:recurp}) by $\rr_i$ on the left gives us
\begin{eqnarray} \label{eq:rpplus1}
\rr_i\pp_{i+1} = -\hf (-2-\ep^2)\pp_i \rr_i + \hf \ep \rr_i\qq_i
\rr_i +\hf \ep (-2-\ep^2)\qq_i +\hf \ep^2 \rr_i\pp_i.
\end{eqnarray}
Rewrite (\ref{eq:recurq}) as
$  \rr_i\qq_i \rr_i = 2q_{i+1}  -  \ep  (p_i \rr_i +\rr_i \pp_i) +
\ep^2 \qq_i.$
Plugging this into (\ref{eq:rpplus1}) and reorganizing the terms,
we obtain (\ref{eq:rpp}). The proof of (\ref{eq:rqq}) is almost
identical.

We will skip the analogous proofs for
(\ref{eq:rpp})--(\ref{eq:rqq}) $\Rightarrow$
(\ref{eq:recurp})--(\ref{eq:recurq}) as well as for
(\ref{eq:recurp})--(\ref{eq:recurq}) $\Rightarrow$
(\ref{eq:rp})--(\ref{eq:rq}).
\end{proof}
%

%
%
\subsection{(Anti-)involutions of $\Haff$}

There are several involutions of the algebra $\Haff$ which are
extensions of the involutions $\sigma, s$ and $-$ of $\Hn$ in
Subsection~\ref{subsec:invol}. We can extend $\sigma$ in two ways to
involutions $\sigma_\pm: \Haff \rightarrow \Haff$ (where $q$ is
fixed):
\begin{eqnarray*}
\sigma_+:& \pp_i \rightarrow \pp_{n+1-i}, &\qq_i \rightarrow
\qq_{n+1-i}, \quad
\rr_i \mapsto \rr_{n-i},  \\
\sigma_-:& \pp_i \rightarrow -\pp_{n+1-i},& \qq_i \rightarrow
-\qq_{n+1-i}, \quad
\rr_i \mapsto \rr_{n-i},
\end{eqnarray*}
for all possible $i$. We extend $s$ to two involutions $s_p, s_q$
of $\Haff$ (where $q$ is fixed):
\begin{eqnarray*}
s_p:& \pp_i \rightarrow -\pp_i, & \qq_i \rightarrow \qq_i, \quad
\rr_i \mapsto -\rr_i, \\
s_q:& \pp_i \rightarrow  \pp_i, & \qq_i \rightarrow -\qq_i, \quad
\rr_i \mapsto -\rr_i,
\end{eqnarray*}
for all possible $i$. We also extend $-$ to  involutions $-_p,
-_q$ of $\Haff$ (fixing each $\rr_i$):
\begin{eqnarray*}
-_p:& \pp_i \rightarrow -\pp_i, & \qq_i \rightarrow \qq_i, \quad
\rr_i \mapsto \rr_i, \quad q \mapsto q^{-1} \\
-_q:& \pp_i \rightarrow  \pp_i, & \qq_i \rightarrow -\qq_i, \quad
\rr_i \mapsto \rr_i, \quad q \mapsto q^{-1},
\end{eqnarray*}
for all possible $i$.
By inspection, all these involutions commute with each other, and
their products give rise to many more involutions of $\Haff$.

Extending the anti-involution $\tau$ on $\Hn$, we also have an
anti-involution $\tau$ on $\Haff$ by letting, for all possible
$i$,
$$\tau (\pp_i) = \pp_i, \quad \tau (\qq_i)
=-\qq_i,\quad \tau (\rr_i) =-\rr_i.
$$
One obtains more anti-involutions on $\Haff$ by composing $\tau$
with the various involutions above (which commute with $\tau$).

\subsection{A nonstandard presentation of the affine Hecke algebra}
The affine Hecke algebra $\affH$ is generated by $T_i $ $(1 \le i
\le n-1), X_j $ $(1\le j \le n)$ subject to the relations
(\ref{hecke})--(\ref{braid}) of $T_i$'s in ${\mathcal H}_n$ and
the following additional relations:
\begin{eqnarray*}
X_iX_j &=& X_jX_i  \;\;\; (\forall i, j)\\
 T_i  X_{i} T_i &=& X_{i+1}\\
 T_i X_j &=& X_j T_i \;\;\; (j \neq i, i+1).
\end{eqnarray*}

Recall $\tprime_i: =T_i +T_i^{-1} \equiv 2T_i -\ep$, and further
introduce
$$P_i := \hf (X_i +X_i^{-1}),
\quad Q_i := \hf (X_i -X_i^{-1}).
$$
It follows that
$$X_i =  P_i + Q_i,
\quad X_i^{-1} = P_i - Q_i.
$$
\begin{proposition}
The algebra $\affH$ is generated by $\tprime_i $ $(1 \le i \le n-1)$, $P_j$ and $Q_j $
$(1\le j \le n)$, subject to the relations
(\ref{eq:tprime})--(\ref{eq:tiplus1prime}) for $\tprime_i$'s and
the following additional relations:
\begin{eqnarray}
\PP_i \PP_j = \PP_j\PP_i,
 && \QQ_i \QQ_j = \QQ_j\QQ_i  \;\; (i \neq j) \label{eq:pq1} \\
\PP_i^2 -\QQ_i^2 =1, && \PP_i \QQ_j = \QQ_j\PP_i, \quad (\forall
i, j)  \label{eq:pq2}
\end{eqnarray}
\begin{eqnarray}
\tprime_i \PP_i &=& \PP_{i+1}\tprime_i -\ep (\QQ_{i+1} +\QQ_i)  \label{eq:pqr1} \\
\tprime_i \QQ_i &=& \QQ_{i+1}\tprime_i -\ep (\PP_i +\PP_{i+1})  \label{eq:pqr2}
\end{eqnarray}
\end{proposition}
\begin{proof}
This follows by a direct computation. Let us illustrate by the derivation of (\ref{eq:pqr1}).
Indeed, recalling  $T_i^{-1} =T_i -\ep$, we have
\begin{eqnarray*}
\tprime_i \PP_i - \PP_{i+1}\tprime_i
&=& \hf(2T_i -\ep) (X_i +X_i^{-1}) - \hf(X_{i+1} +X_{i+1}^{-1})(2T_i -\ep)  \\
&=& T_i X_i -\hf \ep X_i +(T_i -\ep) X_i^{-1} +\hf \ep X_i^{-1}   \\
&&
 -X_{i+1} (T_i -\ep) -\hf \ep X_{i+1} - X_{i+1}^{-1} T_i +\hf \ep X_{i+1}^{-1}   \\
&=& (T_i X_i -X_{i+1} T_i^{-1}) + (T_i^{-1} X_i^{-1} - X_{i+1}^{-1} T_i)
-\ep (\QQ_{i+1} +\QQ_i)  \\
&=& -\ep (\QQ_{i+1} +\QQ_i).
\end{eqnarray*}
In the last equation, we have used $T_i X_i =X_{i+1} T_i^{-1}$,
and $T_i^{-1} X_i^{-1} =X_{i+1}^{-1} T_i.$
\end{proof}

One further checks that in the presence of (\ref{eq:tprime}), the
relations (\ref{eq:pqr1})--(\ref{eq:pqr2}) are equivalent to the two
equations below.
\begin{eqnarray*}
\tprime_i \PP_{i+1} &=&  \PP_{i}\tprime_i +\ep (\QQ_{i+1} +\QQ_i)  \\
\tprime_i \QQ_{i+1} &=&  \QQ_{i}\tprime_i +\ep (\PP_i +\PP_{i+1}).
\end{eqnarray*}

\subsection{The covering affine Hecke algebra}

\begin{definition}
The {\em covering affine Hecke algebra} $\tHaff$ is generated by
$z$, $\tilde{T}_i $ $(1 \le i \le n-1)$, $\tilde{\PP}_j$ and
$\tilde{\QQ}_j $ $(1\le j \le n)$, subject to the relations
(\ref{eq:zt})--(\ref{eq:ztitiplus1}) for $\tilde{T}_i$'s and the
following additional relations:
\begin{eqnarray*}
z^2 =1, && z \text{ is central} \\
\tilde{\PP}_i \tilde{\PP}_j = \tilde{\PP}_j\tilde{\PP}_i,
 && \tilde{\QQ}_i \tilde{\QQ}_j = z\tilde{\QQ}_j\tilde{\QQ}_i  \;\; (i \neq j) \label{eq:zpq1} \\
\tilde{\PP}_i^2 -z\tilde{\QQ}_i^2 =1, && \tilde{\PP}_i
\tilde{\QQ}_j = \tilde{\QQ}_j\tilde{\PP}_i \quad (\forall i, j)
\label{eq:zpq2}
\end{eqnarray*}
\begin{eqnarray}
\tilde{T}_i \tilde{\PP}_i
&=& \tilde{\PP}_{i+1}\tilde{T}_i -\ep (\tilde{\QQ}_{i+1} +z\tilde{\QQ}_i)  \label{eq:zpqr1} \\
\tilde{T}_i \tilde{\QQ}_i &=& z\tilde{\QQ}_{i+1}\tilde{T}_i -\ep
(\tilde{\PP}_i +\tilde{\PP}_{i+1})  \label{eq:zpqr2}
\end{eqnarray}
\end{definition}
Similar to Proposition~\ref{prop:recur}, we can show, by assuming
only (\ref{eq:zt}), that the three pairs of relations
(\ref{eq:zpqr1})--(\ref{eq:zpqr2}),
(\ref{eq:zpqr3})--(\ref{eq:zpqr4}), and
(\ref{eq:zpqr5})--(\ref{eq:zpqr6}) below are equivalent to each
other:
\begin{eqnarray}
\tilde{T}_i \tilde{\PP}_{i+1}
 &=&  \tilde{\PP}_{i}\tilde{T}_i +\ep (\tilde{\QQ}_{i+1} +z\tilde{\QQ}_i)  \label{eq:zpqr3} \\
\tilde{T}_i \tilde{\QQ}_{i+1} &=& z\tilde{\QQ}_{i}\tilde{T}_i +z
\ep (\tilde{\PP}_i +\tilde{\PP}_{i+1}). \label{eq:zpqr4} \\
\tilde{\PP}_{i+1} &=& \frac{1}8 (3 -z)\left(
 z \tilde{T}_i \tilde{\PP}_i \tilde{T}_i
  + \ep  \tilde{T}_i \tilde{\QQ}_i +\ep \tilde{\QQ}_i \tilde{T}_i
  +\ep^2  \tilde{\PP}_i \right)   \label{eq:zpqr5} \\
 \tilde{\QQ}_{i+1} &=& \frac{1}8 (3 -z) \left(
 \tilde{T}_i \tilde{\QQ}_i \tilde{T}_i
  + \ep  \tilde{T}_i \tilde{\PP}_i +\ep \tilde{\PP}_i \tilde{T}_i
  +\ep^2 z \tilde{\QQ}_i \right).   \label{eq:zpqr6}
\end{eqnarray}

By definition, the quotient of the covering affine Hecke algebra
$\tHaff$ by the ideal $\langle z-1 \rangle$ is isomorphic to the
usual affine Hecke algebra in the nonstandard presentation above,
where the canonical images of $\tilde{\PP}_i, \tilde{\QQ}_i$ are
identified with $\PP_i, \QQ_i$ respectively. Also, the quotient of
$\tHaff$ by the ideal $\langle z+1 \rangle$ is isomorphic to the
spin affine Hecke algebra $\Haff$, where the canonical images of
$\tilde{\PP}_i, \tilde{\QQ}_i$ are identified with $\pp_i, \qq_i$
respectively.

\section{Structures of the spin affine Hecke algebra}
\label{sec:structureaffine}
\subsection{Morita super-equivalence for $\Haff$}

\begin{theorem}  \label{th:isomaffine}
There exists an isomorphism of superalgebras
$$
{\Phi}: \HCaff \stackrel{\simeq}{\longrightarrow} \Cl_n \bigotimes
\Haff$$ which extends the isomorphism $\Phi: \HCn \rightarrow
\Cl_n \otimes \Hn$ and is such that
$${\Phi} (X_i) = \pp_i -c_i\qq_i, \quad
{\Phi} (X_i^{-1}) = \pp_i +c_i\qq_i.$$ The inverse ${\Psi}$ is an
extension of $\Psi: \Cl_n \otimes \Hn \rightarrow \HCn$ such that
$${\Psi} (\pp_i) =\hf (X_i +X_i^{-1}), \quad
{\Psi} (\qq_i) =\hf (X_i -X_i^{-1})c_i.
$$
\end{theorem}
\subsection{Proof of the isomorphism Theorem~\ref{th:isomaffine}}
We will adopt the convention ${\Psi} (a) =a^{{\Psi}}, {\Phi} (b)
=b^{{\Phi}}$. We start with several lemmas.
\begin{lemma}  \label{lem:rrppqq}
In $\HCaff$, we have
\begin{eqnarray*}
\rr_i^{{\Psi}} \pp_i^{{\Psi}} &=& \pp_{i+1}^{{\Psi}} \rr_i^{{\Psi}}
 + \ep (\qq_i^{{\Psi}} -\qq_{i+1}^{{\Psi}})
  \label{rrpp} \\
\rr_i^{{\Psi}} \qq_i^{{\Psi}}
 &=& -\qq_{i+1}^{{\Psi}}\rr_i^{{\Psi}}
  -\ep (\pp_i^{{\Psi}} +\pp_{i+1}^{{\Psi}}).
   \label{rrqq}
\end{eqnarray*}
\end{lemma}

\begin{proof}
By (\ref{equivexp}), we have
\begin{eqnarray*}
\rr_i^{{\Psi}} X_i
 &=& (T_i(c_{i+1} -c_i) +\ep c_{i}) X_i \\
 &=& T_i X_i c_{i+1} -T_iX_i^{-1}c_i +\ep X_i^{-1}c_i  \\
 &=& X_{i+1} T_i c_{i+1} -\ep (X_{i+1} +c_ic_{i+1} X_i)c_{i+1} \\
  &&
  -X_{i+1}^{-1} T_i c_i  -\ep (X_i^{-1} c_i -X_{i+1}^{-1}c_{i+1})
   +\ep X_i^{-1} c_i
  \\
 &=& X_{i+1} T_i c_{i+1}  -X_{i+1}^{-1} T_i c_i
  -\ep (X_{i+1} c_{i+1} -X_{i+1}^{-1}c_{i+1} + X_i^{-1} c_i).
\end{eqnarray*}
On the other hand, we have
\begin{eqnarray*}
\rr_i^{{\Psi}} X_i^{-1}
 &=& (T_i(c_{i+1} -c_i) +\ep c_{i}) X_i^{-1} \\
 &=& T_i X_i^{-1} c_{i+1} -T_iX_ic_i +\ep X_ic_i  \\
 &=& X_{i+1}^{-1} T_i c_{i+1} +\ep (X_{i}^{-1} +c_ic_{i+1} X_{i+1})c_{i+1} \\
  &&
  -(X_{i+1} T_i -\ep (X_{i+1}+c_ic_{i+1} X_{i}) )c_{i}
   +\ep X_i c_i
  \\
 &=& X_{i+1}^{-1} T_i c_{i+1}  -X_{i+1} T_i c_i
  +\ep (X_{i+1} c_i +X_{i+1}^{-1}c_i + X_i c_i).
\end{eqnarray*}
Now the lemma follows by adding and subtracting these two
identities for $\rr_i^{{\Psi}} X_i$ and $\rr_i^{{\Psi}} X_i^{-1}$
(as well as  multiplying with $c_i$).
\end{proof}

\begin{lemma} \label{lem:pcqc}
Let $1\le i \le n$. In $\HCaff$, the element $\pp_i^{{\Psi}}$
commutes with $\Cl_n$ while $\qq_i^{{\Psi}}$ super-commutes with
$\Cl_n$.
\end{lemma}
\begin{proof}
Follows directly from (\ref{clifford}) and (\ref{eq:xc}).
\end{proof}

\begin{lemma} \label{lem:txphi}
The $T_i^\Phi, X_i^\Phi, (X_i^-)^\Phi, c_i$ satisfy the relations
(\ref{eq:txt})--(\ref{eq:xc}).
\end{lemma}
\begin{proof}
The relations (\ref{eq:tixj})--(\ref{eq:xc}) for $X_i^\Phi,
(X_i^-)^\Phi, c_i$ are easy to verify from the definitions. It
remains to check (\ref{eq:txt}). We shall use repeatedly
(\ref{eq:rp}) and (\ref{eq:rq}) below. Recalling $T_i^\Phi$ from
(\ref{phiT}), we calculate that
\begin{eqnarray*}
 && 2 X_i^\Phi T_i^\Phi \\
 &=& (\pp_i -c_i\qq_i)
 (\rr_i (c_{i+1}-c_i) +\ep (1 -c_ic_{i+1})) \\
 &=& (\rr_i \pp_{i+1} +\ep (\qq_i -\qq_{i+1})) (c_{i+1} -c_i)
 - (-\rr_i \qq_{i+1} -\ep (\pp_i +\pp_{i+1})) c_i(c_{i+1} -c_i) \\
 &&
 +\ep \pp_i (1-c_ic_{i+1}) +\ep \qq_i c_i (1-c_ic_{i+1})  \\
 &=& (\rr_i \pp_{i+1} -\ep \qq_{i+1}) (c_{i+1} -c_i)
 -(\rr_i\qq_{i+1} +\ep \pp_{i+1}) (1 -c_ic_{i+1}).
\end{eqnarray*}
Therefore, we have
\begin{eqnarray*}
 && 4(T_i^\Phi +\ep c_i c_{i+1}) X_i^\Phi T_i^\Phi \\
 &=& [\rr_i (c_{i+1}-c_i) +\ep (1 -c_ic_{i+1})] \times \\
 &&
  [(\rr_i \pp_{i+1} -\ep \qq_{i+1}) (c_{i+1} -c_i)
 -(\rr_i\qq_{i+1} +\ep \pp_{i+1}) (1 -c_ic_{i+1})]  \\
 &=& -\rr_i (\rr_i \pp_{i+1} -\ep \qq_{i+1}) (c_{i+1} -c_i)^2
  -\rr_i(\rr_i \qq_{i+1}+\ep \pp_{i+1})(c_{i+1}-c_i)(1-c_ic_{i+1})
  \\
  &&
  +\ep (\rr_i\pp_{i+1} -\ep \qq_{i+1})(1+c_ic_{i+1})(c_{i+1} -c_i)
  \\ &&
  -\ep(\rr_i\qq_{i+1} +\ep \pp_{i+1})(1+c_ic_{i+1})(1-c_ic_{i+1}) \\
  &=& 4 (\pp_{i+1} -c_{i+1}\qq_{i+1}) =4 X_{i+1}^\Phi.
\end{eqnarray*}
This completes the proof of the lemma.
\end{proof}

\begin{proof}[Proof of Theorem~\ref{th:isomaffine}]
It is straightforward to check that $\pp_i^\Psi, \qq_i^\Psi,
\rr_i^\Psi$ satisfy (\ref{ppqqcommute})--(\ref{rpqcommute}).
Together with Lemmas~\ref{lem:commute}, \ref{lem:homom},
\ref{lem:braidmod}, \ref{lem:rrppqq}, \ref{lem:pcqc}, this implies
that $\Psi: \Cl_n \otimes \Haff \rightarrow \HCaff$ is an algebra
homomorphism.

Clearly $X_i^\Phi \cdot (X_i^{-1})^\Phi =1$. Recalling that $\Phi
|_{\HCn}: \HCn \rightarrow \Cl_n \otimes \Hn$ is an algebra
isomorphism by Theorem~\ref{th:isomfiniteq}, we have by
Lemma~\ref{lem:txphi} that $\Phi: \HCaff \rightarrow \Cl_n \otimes
\Haff$ is an algebra homomorphism.

By a direct computation, the homomorphisms $\Psi$ and $\Phi$ are
inverses on the generators, and thus they are inverse algebra
isomorphisms.
\end{proof}

\subsection{A basis for $\Haff$}
We recall the definition of $\rr_{i,a_i}$ from (\ref{monomial}).
\begin{theorem}
The algebra $\Haff$ has a basis
$$\pp_1^{k_1} \cdots \pp_n^{k_n} \qq_1^{\epsilon_1} \cdots \qq_n^{\epsilon_n}
\cdot
\rr_{1,a_1}\rr_{2,a_2}\cdots \rr_{n-1,a_{n-1}}
$$
where $k_1,\ldots, k_n \in\Z_+, \epsilon_1,\ldots,\epsilon_n \in
\{0,1\}$, $0 \le a_i\le i$ and $1 \le i \le n-1$.
\end{theorem}

\begin{proof}
The subalgebra $A_i$ generated by $c_i, \pp_i^\Psi, \qq_i^\Psi$
(for a fixed $i$) is identical to the subalgebra generated by
$c_i, X_i, X_i^{-1}$, and it has a linear basis given by
$c_i^\alpha X_i^a$ $(\alpha \in \{0,1\}, a\in \Z)$. By the
standard basis for $\HCaff$, Theorem~\ref{th:isomfiniteq} and
Theorem~\ref{th:isomaffine}, we have the following isomorphisms of
vector spaces:
\begin{eqnarray} \label{eq:twoiso}
\HCaff \cong A_1 \otimes \cdots \otimes A_n \otimes \Hn \cong
\Cl_n \otimes \Haff.
\end{eqnarray}

{\em Claim.} The algebra $A_i$ has another basis $\{c_i^\alpha
(\pp_i^\Psi)^k (\qq_i^\Psi)^\beta\mid \alpha, \beta \in \{0,1\}, k
\in \Z_+\}$. Equivalently, the subalgebra $B_i$ generated by $(X_i
+X_i^{-1})$ and $(X_i -X_i^{-1})$ has a basis $\{(X_i +X_i^{-1})^k
(X_i -X_i^{-1})^\beta \mid \beta \in \{0,1\}, k \in \Z_+\}$.

We prove the second (equivalent) part of the claim. Since any even
power of $(X_i -X_i^{-1})$ can be written as a polynomial in $(X_i
+X_i^{-1})$, the algebra $B_i$ is spanned by the elements $(X_i
+X_i^{-1})^k (X_i -X_i^{-1})^\beta$, with the constraint $\beta \in
\{0,1\}$. It remains to prove the linear independence of these
elements. Assume otherwise
$$f(X):=\sum_k a_k (X_i +X_i^{-1})^k (X_i -X_i^{-1})
 + \sum_k b_k (X_i +X_i^{-1})^k =0,
$$
for some $a_k, b_k \in \C (q)$, all but finitely many  being zero.
Then $f(X^{-1}) =0$, and thus
$$2\sum_k b_k (X_i +X_i^{-1})^k =f(X) +f(X^{-1})=0.$$
By looking at the highest degree term in $X_i$ of this equation, we
see that all $b_k = 0$, and similarly all $a_k = 0$. This proves the
claim.

The theorem now follows from (\ref{eq:twoiso}) and the claim.
\end{proof}
Accordingly, we obtain a similar basis for the covering affine
Hecke algebra $\tHaff$ (compare Proposition~\ref{PBWcovering}).
\subsection{The center of $\Haff$}

\begin{theorem}  \label{th:center}
The even center of the algebra $\Haff$ is the algebra of symmetric
polynomials $\C[p_1, p_2, \cdots, p_n]^{S_n}$.
\end{theorem}

\begin{proof}
The center of $\HCaff$ is equal to $\C[X_1 +X_1^{-1},  \cdots, X_n
+X_n^{-1}]^{S_n}$, the algebra of symmetric polynomials in $X_k
+X_k^{-1} $ $(1\le k \le n)$, according to Jones-Nazarov \cite{JN}
(cf. \cite{BK}). The map $\Phi: \HCaff
\stackrel{\simeq}{\longrightarrow} \Cl_n \otimes \Haff$ sends
$\C[X_1 +X_1^{-1},  \cdots, X_n +X_n^{-1}]^{S_n}$ onto $\C[p_1, p_2,
\cdots, p_n]^{S_n}$, and so $\C[p_1, p_2, \cdots, p_n]^{S_n}$ is the
center of $\Cl_n \otimes \Haff$. It follows that $\C[p_1, p_2,
\cdots, p_n]^{S_n}$ is contained in the (even) center of $\Haff$. On
the other hand, any given even central element $e$ of $\Haff$
commutes with $\Cl_n$ thanks to the evenness of $e$, and thus lies
in the center of $\Cl_n \otimes \Haff$, which is $\C[p_1, p_2,
\cdots, p_n]^{S_n}$.
\end{proof}

\begin{proposition}
Let $n >0$ be odd. Then ${\mathfrak q} :=\qq_1\qq_2\cdots \qq_n$
is an odd central element of $\Haff$. However, ${\mathfrak
q}^\Psi$ does not lie in the center of $\HCaff$.
\end{proposition}

\begin{proof}
${\mathfrak q}^\Psi$ does not lie in the center of $\HCaff$, since
${\mathfrak q}^\Psi c_i =-c_i {\mathfrak q}^\Psi$.

By definition, ${\mathfrak q}$ commutes with each $\pp_i$. Since
$n$ is odd, ${\mathfrak q}$ commutes with each $\qq_i$ by a direct
computation. It remains to show that $\rr_i {\mathfrak q}
={\mathfrak q} \rr_i$ for each $i$. Indeed, by (\ref{eq:rq}) and
(\ref{eq:rqq}), we have
\begin{eqnarray*}
\rr_i \qq_i\qq_{i+1} &=& (-\qq_{i+1} \rr_i -\ep (\pp_i
+\pp_{i+1})) \qq_{i+1}  \\
 &=& -\qq_{i+1} (-\qq_i \rr_i -\ep (\pp_i +\pp_{i+1})) -\ep (\pp_i +\pp_{i+1})
\qq_{i+1} = -\qq_i \qq_{i+1} \rr_i.
\end{eqnarray*}
This, together with (\ref{rpqcommute}) and the oddness of $n$,
implies that $\rr_i {\mathfrak q} ={\mathfrak q} \rr_i$.
\end{proof}

\subsection{The intertwiners}

By Theorem~\ref{th:center}, $\delta :=\prod_{1\le i<j\le n} (\pp_i
-\pp_j)^2$ is an even central element in $\Haff$. Denote by
$(\Haff)_\delta$ the localization of $\Haff$ at $\delta$. In
particular, $(\pp_i -\pp_{i+1})^{-1} \in (\Haff)_\delta$.

 Define
$$\gimel_i =R_i -\ep \frac{\qq_i -\qq_{i+1}}{\pp_i -\pp_{i+1}} \in (\Haff)_\delta.$$
It is understood here and below that $\frac{A}{B} = B^{-1}A$.

\begin{proposition}
 The elements $\gimel_i$ $(1\le i \le n-1)$ satisfy the following
 relations:
\begin{eqnarray}
\gimel_i^2 &=& -2 + 2\ep^2   \frac{\pp_i \pp_{i+1} -1}{(\pp_i
-\pp_{i+1})^2}  \label{gimel2} \\
 \gimel_i \gimel_{i+1} \gimel_i
 &=& \gimel_{i+1} \gimel_i \gimel_{i+1} \label{braidgimel}  \\
\gimel_i\gimel_j &=& - \gimel_j \gimel_i \;\; (|i-j|>1) \nonumber
\end{eqnarray}
\begin{eqnarray*}
\gimel_i \pp_i = \pp_{i+1} \gimel_i, && \gimel_i \pp_{i+1} =
\pp_{i} \gimel_i  \\
\gimel_i \pp_j =\pp_j \gimel_i,   &&
\gimel_i \qq_j = -\qq_j \gimel_i \;\;\, (j \neq i,i+1)   \\
\gimel_i \qq_i = -\qq_{i+1} \gimel_i, && \gimel_i \qq_{i+1} =
-\qq_{i} \gimel_i.
\end{eqnarray*}
\end{proposition}

\begin{proof}
All of these relations can be verified by direct computation. Below
we describe an alternative way by making connections with the
intertwiners $\phi_i$ for $\HCaff$ introduced in \cite[(3.6)]{JN}.
Recall that
 $$\phi_i :=T_i +\frac{\ep}{X_iX_{i+1}^{-1} -1}
 - \frac{\ep}{X_iX_{i+1} -1} \cdot c_ic_{i+1}$$
in a suitable localization of $\HCaff$ isomorphic to $\Cl_n
\otimes (\Haff)_\delta$. One can show that
$$\phi_i c_i = c_{i+1} \phi_i, \quad \phi_i c_{i+1} = c_i \phi_i.$$

{\em Claim.} The isomorphism $\Phi: \HCaff {\longrightarrow} \Cl_n
\otimes \Haff$ sends $\phi_i$ to $\hf (c_i-c_{i+1}) \otimes
\gimel_i$.

Indeed, we have
\begin{eqnarray*}
\Phi (\phi_i) &=&
 \hf \rr_i (c_{i+1} -c_i) +\frac{\ep}{2} (1-c_ic_{i+1}) \\
 &&
  +\ep \Phi \left(
  \frac{X_{i+1} -X_i^{-1} -
  (X_{i+1}^{-1}-X_i^{-1})c_ic_{i+1}}{X_i +X_i^{-1} -(X_{i+1}
  +X_{i+1}^{-1})}
  \right) \\
&=&
 \hf \rr_i (c_{i+1} -c_i)
  +\frac{\ep}{2} \cdot
   \frac{\qq_ic_i +\qq_{i+1}c_{i+1} -\qq_{i+1}c_i -\qq_ic_{i+1}}{\pp_i -\pp_{i+1}}\\
&=&
 \hf \left (\rr_i - \ep \frac{\qq_i-\qq_{i+1}}{\pp_i -\pp_{i+1}} \right )
   (c_{i+1} -c_i)
   = \hf (c_i-c_{i+1}) \otimes \gimel_i.
\end{eqnarray*}

With the help of the claim, all of the identities in the Proposition
follow from the corresponding statements for $\phi_i$ in
\cite[(3.7)]{JN} and \cite[Prop.~3.1]{JN}. Let us illustrate by
proving (\ref{gimel2}) in detail below. Recall from
\cite[Prop.~3.1]{JN} that
\begin{eqnarray*}
\phi_i^2 = 1 -\ep^2 \left (
\frac{X_iX_{i+1}^{-1}}{(X_iX_{i+1}^{-1} -1)^2}
 + \frac{X_i^{-1}X_{i+1}^{-1}}{(X_i^{-1}X_{i+1}^{-1} -1)^2}
\right ).
\end{eqnarray*}
By the above Claim, we have
\begin{eqnarray*}
\gimel_i^2 &=& \left( (c_i-c_{i+1}) \Phi (\phi_i) \right)^2
= -2 \Phi (\phi_i^2) \\
 &=& -2  +2 \ep^2
 \frac{X_i^{-1}X_{i+1}^{-1} -2 +X_iX_{i+1} +X_iX_{i+1}^{-1} -2
+X_i^{-1}X_{i+1}}{(X_{i+1}+X_{i+1}^{-1} -X_i-X_i^{-1})^2} \\
 &=& -2  +2 \ep^2
 \frac{(X_{i+1}+X_{i+1}^{-1})
(X_i+X_i^{-1}) -4}{(X_{i+1}+X_{i+1}^{-1} -X_i-X_i^{-1})^2} \\
 &=&
-2 + 2\ep^2 \frac{\pp_i \pp_{i+1} -1}{(\pp_i -\pp_{i+1})^2}.
\end{eqnarray*}

For the braid relation (\ref{braidgimel}), the following identity
can be useful:
$$(c_i -c_{i+1})(c_{i+1}-c_{i+2})(c_i -c_{i+1})
 =(c_{i+1}-c_{i+2})(c_i -c_{i+1})(c_{i+1}-c_{i+2}).$$
\end{proof}

\section{Cyclotomic spin Hecke algebras}
\label{sec:cyclot}

\subsection{The definition}

Recall $\pp_1\qq_1 =\qq_1\pp_1$. Consider the subalgebra
$$\mathcal A_1 :=\C[\pp_1, \qq_1]/\langle \pp_1^2 +\qq_1^2-1
\rangle$$
of $\Haff$ which is commutative and $\Z_2$-graded with $\pp_1$
being even and $\qq_1$ odd.

\begin{proposition} \label{prop:ideal}
A nonzero $\Z_2$-homogeneous ideal $I_1$ of $\mathcal A_1$ is one
of the following:
\begin{enumerate}
\item
  $\langle f(\pp_1) \rangle$, for some nonzero polynomial $f$ in one variable;
\item
  $\langle g(p_1)\qq_1 \rangle$, for some nonzero polynomial $g$ in one variable;
\item
 $\langle (\pp_1 + 1) g(\pp_1), g(\pp_1) \qq_1 \rangle$, for some nonzero polynomial $g$;
\item
 $\langle (\pp_1 - 1) g(\pp_1), g(\pp_1) \qq_1 \rangle$, for some nonzero polynomial
$g$.
\end{enumerate}
\end{proposition}

\begin{proof}

Let $I_1$ be a nonzero $\Z_2$-homogeneous ideal of $\mathcal A_1$.
Let $f$ and $g$ be the unique monic polynomials of minimal degree
such that $f(\pp_1) \in I_1, g(\pp_1) \qq_1 \in I_1$. By the
$\Z_2$-homogeneity, $I_1 =\langle f(\pp_1), g(\pp_1) \qq_1 \rangle.$

Note that $f(\pp_1) \qq_1 \in I_1,$ and $(\pp_1^2-1) \cdot
g(\pp_1) =-g(\pp_1) \qq_1^2 \in I_1$. By assumption of minimal
degrees on $f,g$, we have
\begin{eqnarray} \label{divisable}
f(\pp_1) \mid (\pp_1-1)(\pp_1+1) \cdot g(\pp_1),
\end{eqnarray}
and thus $\deg f \le \deg g+2$. Also $\deg g \le \deg f$, and
$g=f$ if $\deg g = \deg f$.

In the case when $\deg f =\deg g$ and thus $g=f$, the ideal $I_1$
is of the form~(1).

In the case when $\deg f =\deg g+2$, we have $f(\pp_1) =
(\pp_1^2-1) \cdot g(\pp_1)$ by (\ref{divisable}), and thus $I_1$
is of the form~(2).

Finally assume that $\deg f =\deg g+1$ and consider two subcases:
(i) $g(\pp_1) \mid f(\pp_1)$; (ii) $g(\pp_1) \nmid f(\pp_1)$. Thanks
to (\ref{divisable}), in case~(i), $I_1$ is of the form~(3) or (4).
We now claim the subcase~(ii) is empty. Indeed, by
(\ref{divisable}), (ii) and $\deg f =\deg g+1$, we have $f(\pp_1)
=(\pp_1^2 -1)h(\pp_1) =-\qq_1^2 h(\pp_1)$ and $g(\pp_1) =(\pp_1 -a)
h(\pp_1)$ for some constant $a \neq \pm 1$ and some polynomial $h$
of degree equal to $(\deg g-1)$. Therefore,
$$ h(\pp_1) \qq_1
 =\frac1{1-a^2}\left( (\pp_1+a) \cdot g(\pp_1)\qq_1 -\qq_1 \cdot f(\pp_1)
\right) \in I_1.$$
This contradicts with the choice of $g(\pp)$ of minimal degree.
\end{proof}

\begin{definition} \label{def:cycspin}
The {\em cyclotomic spin Hecke algebra} $\Hcycl$ is the quotient
algebra of $\Haff$ by the two-sided ideal $I =\langle I_1 \rangle$
generated by a nonzero $\Z_2$-homogeneous ideal $I_1 \subset
\mathcal A_1$. (Note that $\Hcycl$ inherits a superalgebra
structure from $\Haff$.)
\end{definition}

\begin{remark}  \label{rem:ideal}
As a byproduct of the above proof of Proposition~\ref{prop:ideal},
the ideal $I$ in Definition~\ref{def:cycspin} is generated by
$f(\pp_1)$ and $g(\pp_1) \qq_1$, where $f$ and $g$ are the unique
monic polynomials of minimal degree such that $f(\pp_1) \in I_1,
g(\pp_1) \qq_1 \in I_1$. We will sometimes write $f=f_I$ and
$g=g_I$ to indicate its dependence on $I$. More specifically, $I$
is generated by one or two elements given in
Proposition~\ref{prop:ideal}.
\end{remark}

\subsection{Relation to cyclotomic Hecke-Clifford
algebras}

We refer to Ariki \cite{Ar} for more on the classical cyclotomic
Hecke algebras.

Let $F$ be a polynomial of the form
$$F(X_1) =a_d X_1^d +a_{d_1}X_1^{d-1} +\cdots a_1X_1+a_0$$
which satisfies the condition
\begin{eqnarray} \label{technical}
a_d =1, \quad a_i =a_0 a_{d-i} \;\; (\forall 0\le i \le d).
\end{eqnarray}
Associated to such an $F$, Brundan-Kleshchev \cite{BK} introduced
the cyclotomic Hecke-Clifford algebra, which will be denoted by
$\HCn^F$ in this paper, as the quotient algebra $\HCaff/\langle
F(X_1)\rangle$. The technical condition (\ref{technical}) was
imposed so that the resulting cyclotomic algebra $\mathcal H_n^F$
has an expected basis and dimension. The next proposition shows
that the condition (\ref{technical}) is natural from the
perspective of cyclotomic spin Hecke algebras.

\begin{theorem} \label{cyclot}
There is a bijection between the set of cyclotomic spin Hecke
algebras and the set of cyclotomic Hecke-Clifford algebras. More
explicitly, every cyclotomic Hecke-Clifford algebra $\HCn^F$ is
isomorphic to $\Cl_n \otimes \Hcycl$ for some cyclotomic spin
Hecke algebra $\Hcycl$ via $\Phi$. Conversely, for each $\Hcycl$,
the algebra $\Cl_n \otimes \Hcycl$ is isomorphic to some
cyclotomic Hecke-Clifford algebra via $\Psi$.
\end{theorem}

\begin{proof}
Note that $a_0 =\pm 1$ by (\ref{technical}). Divide the degree $d$
polynomials $F$ which satisfy the condition~(\ref{technical}) into
the following four cases:
\begin{enumerate}
 \item
$d=2k$ is even and $a_0 =1$;
 \item
$d=2k$ is even and
$a_0 =-1$;
 \item
$d=2k+1$ is odd and $a_0=1$;
 \item $d=2k+1$ is odd and $a_0= -1$.
\end{enumerate}
Then it follows by a case-by-case elementary verification that the
isomorphism $\Phi: \HCaff \rightarrow \Cl_n\otimes \Haff$ sends
$X_1^{-k} F(X_1)$ for $F$ in each case bijectively onto the
corresponding set below:
\begin{enumerate}
 \item
 $\{f(\pp_1) \mid f \text{ is a polynomial of degree } k\}$;
 \item
 $\{g(\pp_1)\qq_1 \mid g \text{ is a polynomial of degree } (k-1)\}$;
 \item
 $\{(\pp_1 +1-c_1\qq_1)\cdot g(\pp_1) \mid g \text{ is a polynomial of degree } k\}$;
 \item
 $\{(\pp_1 -1-c_1\qq_1)\cdot g(\pp_1) \mid g \text{ is a polynomial of degree }
k\}$.
\end{enumerate}

Clearly, the ideal $\mathfrak I$ in $\Cl_n\otimes \Haff$ generated
by an element in (1) or (2) above coincides with $\Cl_n \otimes
\langle I_1\rangle$ where $I_1$ is  given by
Prop.~\ref{prop:ideal}~(1) or (2) respectively. Now the proposition
follows by the following claim.

{\em Claim.} The ideal $\mathfrak I$ in $\Cl_n\otimes \Haff$
generated by the element $(\pp_1\pm 1-c_1\qq_1)\cdot g(\pp_1)$ in
(3) or (4) coincides with $\Cl_n \otimes \langle I_1\rangle$, where
$\langle I_1\rangle$ is the ideal in $\Haff$ generated by $I_1$ in
Proposition~\ref{prop:ideal}~(3) or (4) respectively.

Let us prove the claim for (3) and skip a similar proof for (4).
Indeed, it is clear that $\mathfrak I \subseteq \Cl_n \otimes
\langle I_1\rangle$. On the other hand, we have
$$ (\pp_1+1)g(\pp_1) = \hf(\pp_1 +1-c_1\qq_1)(\pp_1
+1+c_1\qq_1)g(\pp_1) \in \mathfrak I, $$
and thus also $ g(\pp_1) \qq_1 =c_1(\pp_1 +1+c_1\qq_1)g(\pp_1)
-c_1(\pp_1+1)g(\pp_1) \in \mathfrak I$. Therefore, $\mathfrak I
\supseteq \Cl_n \otimes \langle I_1\rangle$.
\end{proof}

It is known \cite{BK} that $\dim \HCn^F = (\deg F)^n 2^nn!.$ From
the explicit relations between (the generators of) the
corresponding ideals in $\HCaff$ and $\Haff$ presented in the
above proof, we have the following.
\begin{corollary}
Let $f_I$ and $g_I$ be the unique monic polynomials of minimal
degree such that $f_I (\pp_1)$ and $g_I (\pp_1) \qq_1$ generate
the ideal $I$ in $\Haff$. Then, $\dim \Hcycl =   (\deg f_I +\deg
g_I)^n n!.$
\end{corollary}

Conjecturally, a basis for $\Hcycl$ consists of
$ \pp_1^{\alpha_1} \qq_1^{\epsilon_1} \cdots \pp_n^{\alpha_n}
\qq_n^{\epsilon_n} \rr$,
where $\epsilon_1, \ldots, \epsilon_n \in \{0,1\}$, $0\leq
\alpha_i <\deg f_I$ if $\epsilon_i=0$ and $0\leq \alpha_i <\deg
g_I$ if $\epsilon_i=1,$ and $\rr$ runs over all standard monomials
in $\Hn$.

%

%
%
\subsection{Jucys-Murphy elements for $\Hn$}

We observe that the spin Hecke algebra $\Hn$ coincides with the
(smallest) cyclotomic spin Hecke algebra $\Hcycl$, where $I
=\langle \pp_1 -1, \qq_1 \rangle$. Similarly, the Hecke-Clifford
algebra $\HCn$ is a special case of the cyclotomic Hecke-Clifford
algebras $\HCn^F$ with $F(X_1) =X_1-1$.

\begin{proposition} \label{prop:JM}
There exists a unique algebra homomorphism
$$\mathcal{JM}: \Haff \longrightarrow \Hn$$
which extends the identity map on $\Hn$ and is such that
$\mathcal{JM} (p_1) =1$, $\mathcal{JM} (q_1) =0$.
\end{proposition}
\begin{proof}
There exists a unique algebra homomorphism $\text{JM}: \HCaff
\rightarrow \HCn$, which extends the identity map on $\HCn$ and is
such that $\text{JM} (X_1) =1$, according to Jones-Nazarov
\cite[Proposition~3.5]{JN}. By (\ref{eq:txt}), the images $J_i$ of
$X_i $ $(1\le i \le n)$ under $\text{JM}$, called the Jucys-Murphy
elements for $\HCaff$, are given recursively by
$J_{i+1} =(T_i +\ep c_ic_{i+1}) J_i T_i$. By
Theorems~\ref{th:isomfiniteq} and \ref{th:isomaffine}, there
exists a homomorphism $\text{JM}': \Haff \rightarrow \Hn$ to make
the following diagram commutative:
$$\CD \HCaff @>\text{JM}>> \HCn \\
@V\Phi V\cong V @V\Phi V\cong V \\
\Cl_n \bigotimes \Haff @> \text{JM}'>> \Cl_n \bigotimes \Hn \endCD
$$

Since $\text{JM} (X_1) =1$, it follows by definition of $\Phi$
that $\text{JM}' (p_1) =1$, $\text{JM}' (q_1) =0$. Moreover, since
$\text{JM}'|_{\Cl_n \otimes \Hn}$ is the identity and the images
of $p_i, q_i$ are given recursively by
Proposition~\ref{prop:recur}, we conclude that $\text{JM}'$ is of
the form $\text{I} \otimes \mathcal{JM}$ for a unique homomorphism
$\mathcal{JM}: \Haff \rightarrow \Hn$ with given images of $\pp_1$
and $\qq_1$.  Note that $\mathcal{JM} (p_1) =1$ and $\mathcal{JM}
(q_1) =0$.
\end{proof}

We will call the images $\mathfrak p_i, \mathfrak q_i \in \Hn$
$(1\le i\le n)$ of the elements $p_i, q_i$'s under the homomorphism
$\mathcal{JM}$ the {\em Jucys-Murphy elements} for $\Hn$, following
the convention for the symmetric group and the usual Hecke algebras.
The relations (\ref{ppqqcommute})--(\ref{eq:rqq}), with $\mathfrak
p_i$ and $\mathfrak q_i$ replacing $p_i$ and $q_i$, are satisfied.
Alternatively, it follows from the proof of
Proposition~\ref{prop:JM} that
$$\mathfrak p_i =\hf \Phi(J_i +J_i^{-1}),
\quad \mathfrak q_i  =\hf \Phi((J_i -J_i^{-1})c_i).
$$
Note the nontrivial implication that $\Phi(J_i +J_i^{-1})$ and
$\hf \Phi((J_i -J_i^{-1})c_i)$ lie in $\Hn$. A direct computation
using the recursive formula in Proposition~\ref{prop:recur} gives
us the first few cases of the Jucys-Murphy elements:
\begin{eqnarray*}
1 =\mathfrak p_1, && \mathfrak q_1 =0, \\
1 +\ep^2 =\mathfrak p_2, &&\mathfrak q_2 =\ep \rr_1, \\
\frac{\ep^2}2 (\rr_1 \rr_2 +\rr_2 \rr_1) + (1+\ep^2)^2 = \mathfrak
p_3, &&
 \mathfrak q_3 =\frac{\ep}2 \left(\rr_1 \rr_2 \rr_1 +
 (2+\ep^2) \rr_2 \right).
\end{eqnarray*}
%
These elements will play important roles in analyzing further the
structures and the representation theory  of $\Hn$ as in the usual
(non-spin) setup.

\subsection{A degeneration of $\Haff$ and $\Hcycl$}
\label{subsec:degen}

Recall that the spin symmetric group algebra $\C S^-_n$ is
generated by $t_i $ $(1\le i \le n-1)$ subject to the relations
(\ref{eq:braidt})--(\ref{eq:tcomm}). The degenerate spin affine
Hecke algebra $\widehat{\mathcal B}$, introduced in \cite{Wa}, is
the superalgebra with odd generators ${b}_i $ $(1\le i \le n)$ and
$t_i $ $(1\le i \le n-1)$, subject to the relations
(\ref{eq:braidt})--(\ref{eq:tcomm}) for $t_i$'s and the following
additional relations:
\begin{eqnarray*}
{b}_i {b}_j &=& - {b}_j {b}_i \quad (i \neq j)  \\
  t_i {b}_i &=& - {b}_{i+1} t_i +1   \\
t_i {b}_j &=& - {b}_j t_i   \quad (j \neq i, i+1).
\end{eqnarray*}

\begin{remark}  \label{rem:deg}
The algebra $\widehat{\mathcal B}$ can be obtained from $\Haff$ by
a suitable degeneration. Set $q =e^{\hbar/2}$. As $q$ goes to $1$,
keeping in mind $\pp_i^2 +\qq_i^2=1$, we set
$$\pp_i \approx 1 + \hbar^2 b_i^2
+o(\hbar^2), \quad \qq_i \approx \hbar \sqrt{-2}\cdot b_i
+o(\hbar), \quad \rr_i \approx \sqrt{-2}\cdot t_i +o(\hbar).$$
Then, as $q$ goes to $1$, the defining relations
(\ref{eq:rri2})--(\ref{braidspin}),
(\ref{ppqqcommute})--(\ref{rpqcommute}), (\ref{eq:rq}) for $\Haff$
reduce to the defining relations for $\widehat{\mathcal B}$. The
remaining relation (\ref{eq:rp}) for $\Haff$ reduces to
%
$ t_i b_i^2 = b_{i+1}^2 t_i +  (b_i -b_{i+1}), $
which follows from the defining relations for $\widehat{\mathcal
B}$.
\end{remark}

\begin{remark}
The isomorphism in Theorem~\ref{th:isomaffine} degenerates in the
sense of Remark~\ref{rem:deg} to the superalgebra isomorphism
between the degenerate affine Hecke-Clifford algebra and  $\Cl_n
\otimes \widehat{\mathcal B}$ established in \cite{Wa}.
\end{remark}

We define the {\em degenerate cyclotomic spin Hecke algebras} as
the quotient algebras ${\mathcal B}^f :=\widehat{\mathcal
B}/\langle f(b_1) \rangle$, where $f$ is an even or an odd
polynomial in one variable. The condition on $f$ is precisely such
that ${\mathcal B}^f$ inherits a canonical superalgebra structure
from $\widehat{\mathcal B}$. Using the Morita super-equivalence
\cite{Wa} between $\widehat{\mathcal B}$ and Nazarov's degenerate
affine Hecke-Clifford algebra, it is straightforward to see that
the degenerate cyclotomic spin Hecke algebras correspond
bijectively to the degenerate cyclotomic Hecke-Clifford algebras
\cite{BK, Kle} (called the cyclotomic Sergeev algebras in {\em
loc. cit.}) via a Morita super-equivalence (compare
Theorem~\ref{cyclot}).

\end{document}